\documentclass[final]{siamltex1213}
\usepackage{layout}
\usepackage{amsmath,amsfonts,amssymb}
\usepackage{xcolor}
\usepackage{hyperref}
\usepackage{epsfig}
\usepackage{todonotes}
\usepackage{algorithm}
\usepackage{algpseudocode}
\usepackage{lscape}
\usepackage{subcaption}
\usepackage{multirow}

\usepackage{pgfplots}
\usepackage[font=sl,labelfont=normal,size=small]{caption}
\usepackage{rotating}

\title{Analyzing the effect of local rounding error propagation on the maximal attainable accuracy of the pipelined Conjugate Gradient method}
\author{
Siegfried Cools\thanks{Applied Mathematics Group, Department of Mathematics and Computer Science, University of Antwerp, Middelheimlaan 1, 2020 Antwerp, BE} 
\and Emrullah Fatih Yetkin\thanks{Department of Computer Engineering, Istanbul Kemerburgaz University, 34217 Ba\u{g}c{\i}lar, TR}
\and Emmanuel Agullo\thanks{HiePACS, Inria Bordeaux - Sud-Ouest, 200 Avenue de la Vieille Tour, 33405 Talence, FR}
\and Luc Giraud\footnotemark[3]
\and Wim Vanroose\footnotemark[1]
} 

\begin{document}

\maketitle

\begin{abstract}
  Pipelined Krylov subspace methods typically offer improved strong scaling on parallel HPC hardware
  compared to standard Krylov subspace methods for large and sparse linear systems. 
	In pipelined methods the traditional synchronization bottleneck is mitigated 
  by overlapping time-consuming global communications with useful computations. 
	However, to achieve this communication hiding strategy, pipelined methods
  introduce additional recurrence relations for a number of auxiliary variables that are required to update the approximate solution. 
	This paper aims at studying the influence of local rounding errors 
	that are introduced by the additional recurrences in the pipelined Conjugate Gradient method.
	Specifically, we analyze the impact of local round-off effects on the attainable accuracy of the pipelined CG algorithm 
	and compare to the traditional CG method. Furthermore, we estimate the gap between the true
	residual and the recursively computed residual used in the algorithm.
	Based on this estimate we suggest an automated residual 
  replacement strategy to reduce the loss of attainable accuracy on the final iterative solution.
  The resulting pipelined CG method with residual replacement improves the maximal 
  attainable accuracy of pipelined CG, while maintaining the efficient parallel performance of the pipelined method.
	This conclusion is substantiated by numerical results for a variety of benchmark problems.
\end{abstract}

\begin{keywords} Conjugate gradients, Parallelization, Latency hiding, Global communication, Pipelining, Rounding errors, Maximal attainable accuracy
\end{keywords}

\section{Introduction} \label{sec:introduction}

Krylov subspace methods \cite{greenbaum1997iterative,liesen2012krylov,meurant1999computer,saad2003iterative,van2003iterative} form the basis linear algebra solvers for many contemporary high-performance computing applications. The Conjugate Gradient (CG) method, which dates back to the 1952 paper by Hestenes and Stiefel \cite{hestenes1952methods},
can be considered as the first of these methods. Although over 60 years old, the CG method is still the work horse method for 
the solution of linear systems with symmetric positive definite (SPD) matrices due to its numerical simplicity and easy implementation. 
These SPD systems may originate from various applications such as e.g.\,the discretization of partial differential equations (PDEs).

Due to the transition of hardware towards the exascale regime in the coming years, research on the scalability of Krylov subspace methods
on massively parallel architectures has become increasingly prominent \cite{dongarra2011international,dongarra1998numerical}. 
This is reflected in the new High Performance Conjugate Gradients (HPCG) benchmark for ranking 
HPC systems introduced by Dongarra et al.~in 2013 \cite{dongarra2013toward, dongarra2015hpcg}. The ranking is based on sparse matrix-vector 
computations and data access patterns, rather than the dense matrix algebra used in the traditional High Performance LINPACK (HPL) benchmark.
Moreover, since the system matrix is often sparse, the main bottleneck for efficient parallel execution is typically not the sparse 
matrix-vector product (\textsc{spmv}), but the communication overhead (bandwidth saturation) caused by global reductions required 
in the computation of dot-products. 

Recently significant research has been devoted to the mitigation and/or elimination of the synchronization bottleneck in Krylov subspace methods.
The earliest papers on synchronization reduction and latency hiding in Krylov subspace methods date back to the late 1980's \cite{strakovs1987effectivity} 
and early 1990's \cite{barrett1994templates,d1992reducing,de1991parallel,demmel1993parallel,erhel1995parallel}.
The idea of reducing the number of global communication points in Krylov subspace methods on parallel computer 
architectures was also used in the $s$-step methods by Chronopoulos 
et al.~\cite{chronopoulos1989s,chronopoulos2010block,chronopoulos1996parallel} and more recently by Carson et al.~in \cite{carson2014residual,carson2013avoiding}. 
In addition to communication avoiding methods\footnote{Although commonly used in the contemporary literature, the term `communication-avoiding' Krylov subspace algorithm it slightly dubious, since the number of global synchronization phases is in fact \emph{reduced} by reformulating the algorithms, rather than \emph{avoided}; hence, the term `communication-reducing' algorithm may be more appropriate in this context.}, research on hiding global communication by overlapping communication with computations was performed by a various authors over the last decades, see Demmel et al.~\cite{demmel1993parallel}, De Sturler et al.~\cite{de1995reducing}, and Ghysels et al.~\cite{ghysels2013hiding,ghysels2014hiding}. We refer the reader to the recent work \cite{carson2016numerical}, Section 2 and the references therein for more background and a wider historical perspective on the development of early variants of the CG algorithm that contributed to the current algorithmic strive towards parallel efficiency.

The pipelined CG (p-CG) method proposed in \cite{ghysels2014hiding} aims at hiding the global synchronization latency of standard 
preconditioned CG by removing some of the global synchronization points. Pipelined CG performs only one global reduction 
per iteration. Furthermore, this global communication phase is overlapped by the sparse 
matrix-vector product (\textsc{spmv}), which requires only local communication. In this way, idle core time 
is minimized by performing useful computations simultaneously to the time-consuming global communication phase, cf.~\cite{eller2015non}.
 
The reorganization of the CG algorithm that is performed to achieve the overlap of communication with computations
introduces several additional  \textsc{axpy} ($y \leftarrow \alpha x +y$) operations to recursively compute auxiliary variables. 
Vector operations such as an \textsc{axpy} are typically computed locally, and thus do not require communication between nodes. 
Thus, the addition of extra recurrences has no impact on the communication flow of the algorithm.
Dot-products of two vectors, on the other hand, involve global communication between all processes, and are therefore grouped together 
in p-CG. 

In exact arithmetic, the resulting pipelined CG algorithm is equivalent to classical CG. However, when switching to finite precision, 
each of the additional recurrences introduce local rounding errors. The propagation of these rounding errors throughout the algorithm
is much more pronounced for pipelined CG, and can have a detrimental effect on the iterative solution.
As a result, a significant loss of attainable accuracy 
compared to classical CG can in practice be observed for the p-CG algorithm.
The current paper contributes to the analysis of the rounding error propagation in different variants of the CG algorithm.
Additionally, the analytical results will be used to formulate an automated residual replacement strategy that 
improves the maximal attainable accuracy of the pipelined CG method. We stress that the proposed residual replacement strategy 
only accounts for an improvement of the attainable accuracy in pipelined CG. Other notable rounding error effects in multi-term recurrence algorithms, such as a delay of convergence due to loss of orthogonality in finite precision arithmetic, see \cite{carson2016numerical} and \cite[Chapter 5]{liesen2012krylov}, are not resolved by the methodology proposed in this work. 

The paper is structured as follows. In Section \ref{sec:analysis} the propagation of local rounding errors in standard preconditioned CG, 
Chronopoulos/Gear CG, and the pipelined CG algorithm is analyzed. Bounds for the gap between the explicitly computed residual and the
recursive residual are derived. Section \ref{sec:practical} proposes an approximate and practically useable estimate for the residual gap.
Furthermore, the incorporation of a residual replacement strategy in the pipelined CG method is discussed in Section \ref{sec:pipelined}. 
A criterion for automated residual replacement based on the aforementioned error propagation model is suggested. 
Extensive numerical experiments in Section \ref{sec:numerical} illustrate the error estimate and show the possible improvement in  
attainable accuracy for the pipelined CG method with automated residual replacement. 
Parallel scaling results indicate that the residual replacement strategy does not affect the parallel scalability of the pipelined CG method. 
Finally, conclusions are formulated in Section \ref{sec:conclusions}. 

\section{Analysis of local rounding error propagation in variants of the CG algorithm} \label{sec:analysis}

The analysis in this section is based upon the rounding error analysis performed by Greenbaum in \cite{greenbaum1997estimating} and Strako{\v{s}} \& Gutknecht in \cite{gutknecht2000accuracy}. Additional work on this topic can be found in \cite{demmel1997applied,greenbaum1989behavior,meurant2006lanczos,paige1971computation,paige1972computational,paige1976error,paige1980accuracy,strakovs2002error,strakovs2005error}. We assume the following classical model for floating point arithmetic on a machine with machine precision $\epsilon$:
\begin{equation}
	\text{fl}(a\pm b) = a(1+\epsilon_1) \pm b(1+\epsilon_2), \quad |\epsilon_1|,|\epsilon_2| \leq \epsilon,
\end{equation}
\begin{equation}
	\text{fl}(a \text{~op~} b) = (a \text{~op~} b) (1+\epsilon_3), \quad |\epsilon_3| \leq \epsilon, \quad \text{op} = *, /,
\end{equation}
where $\text{fl}(a)$ denotes the finite precision floating point representation of the mathematical quantity $a$. Under this model, and discarding terms involving $\epsilon^2$ or higher powers of $\epsilon$ when terms of order $\epsilon$ are present, the following standard results for operations on an $n$-by-$n$ matrix $A$, $n$-length vectors $v$ and $w$ and scalar number $\alpha$ hold:
\begin{equation}
	\| \alpha v - \text{fl}(\alpha v) \| \leq \| \alpha v \| \, \epsilon =  |\alpha| \, \|v\| \, \epsilon,
\end{equation}
\begin{equation}
	\| v + w - \text{fl}(v + w) \| \leq (\|v\| + \|w\|) \, \epsilon,
\end{equation}
\begin{equation}
	| \left( v,w \right) - \text{fl}(\,\left(v,w \right)\,) | \leq n \, \|v\| \, \|w\| \epsilon,
\end{equation}
\begin{equation}
	\| Av - \text{fl}(Av) \| \leq (\mu\sqrt{n}) \, \|A\| \, \|v\| \, \epsilon ,
\end{equation}
where $\mu$ is the maximum number of nonzeros in any row of $A$. The norm $\|\cdot \|$ denotes the Euclidean 2-norm throughout this manuscript, unless explicitly stated otherwise.

\subsection{Accumulation of local rounding errors in classical CG}

Classical preconditioned CG is given by Algorithm \ref{algo::pcg}. Note that in the unpreconditioned case, line 8 is dropped, and each occurrence of $u_i$ is replaced by $r_i$.
In finite precision arithmetic, the recurrences for the computed search direction $\bar{p}_i$, iterate $\bar{x}_i$ and residual $\bar{r}_i$ in iteration $i$ ($i = 0,1,2,\ldots$) of the CG algorithm are
\begin{eqnarray} 
	\bar{p}_{i+1} &=& \bar{u}_{i+1} + \bar{\beta}_{i+1} \bar{p}_i + \delta_i^p , \notag \\
	\bar{x}_{i+1} &=& \bar{x}_i + \bar{\alpha}_i \bar{p}_i + \delta_i^x , \notag \\
	\bar{r}_{i+1} &=& \bar{r}_i - \bar{\alpha}_i A \bar{p}_i + \delta_i^r \label{eq:pxr},
\end{eqnarray}
where $\delta_i^p$, $\delta_i^r$ and $\delta_i^x$ contain the local rounding errors produced in step $i$. In our notation barred variables (e.g., $\bar{p}_i$, $\bar{x}_i$, $\bar{r}_i$, $\bar{u}_i$, $\bar{\alpha}_i$ and $\bar{\beta}_{i+1}$) will always denote the actually computed quantities. This notation should avoid confusion with the mathematical quantities defined in exact arithmetic (e.g., $r_i$, $p_i$ and $\alpha_i$) which are unavailable in practice. Vectors obtained by actually applying the matrix-vector product will be referred to as \emph{explicit} quantities (e.g., $b-A\bar{x}_i$ is called the \emph{explicit} or \emph{true residual}), in contrast to the \emph{recursive} quantities given by \eqref{eq:pxr} (e.g., $\bar{r}_i$ is called the \emph{recursive residual}). Since the \textsc{spmv} is a computationally costly operation, the residual is only computed recursively in Algorithm \ref{algo::pcg}, except for $\bar{r}_0 = \text{fl}(b-A\bar{x}_0)$.

\begin{algorithm}[t]
  \caption{Preconditioned CG}
  \label{algo::pcg}
  \begin{algorithmic}[1]
  	\Procedure{prec-cg}{$A$, $M^{-1}$, $b$, $x_0$}
    \State $r_0 := b - Ax_0$; $u_0:= M^{-1} r_0$ ; $p_0 = u_0$  \label{al:residual}
    \For{$i = 0, \dots$}
    \State $s_i := Ap_{i}$ \label{al:matvec}
    \State $\alpha_{i} := \left( r_i, u_i \right) / \left( s_i, p_i \right)$ \label{al:alpha}
    \State $x_{i+1} := x_i + \alpha_{i} p_i$
    \State $r_{i+1} := r_i - \alpha_{i} s_i$
    \State $u_{i+1} := M^{-1} r_{i+1}$
    \State $\beta_{i+1} := \left( r_{i+1}, u_{i+1} \right) / \left( r_i, u_i \right)$
    \State $p_{i+1} := u_{i+1} + \beta_{i+1} p_i$ 
    \EndFor
    \EndProcedure
  \end{algorithmic}
\end{algorithm}

The iteration $i$ local rounding errors satisfy the following bounds:
\begin{eqnarray} 
	\|\delta_i^p\| &\leq& \left( \|\bar{u}_{i+1}\| + 2 \, |\bar{\beta}_{i+1}| \, \|\bar{p}_i\| \right) \epsilon, \notag \\
	\|\delta_i^x\| &\leq& \left( \|\bar{x}_i\| + 2 \, |\bar{\alpha}_i| \, \|\bar{p}_i\| \right) \epsilon , \notag \\
	\|\delta_i^r\| &\leq& \left( \|\bar{r}_i\| + (\mu\sqrt{n}+2) \, |\bar{\alpha}_i| \, \|A\| \, \|\bar{p}_i\| \right) \epsilon \label{eq:delta_pxr}.
\end{eqnarray}
We now want to estimate the difference (or gap) between the true residual $b-A\bar{x}_i$ and the recursive residual $\bar{r}_i$. Hence, we define this gap as 
\begin{equation}
	f_i = (b-A\bar{x}_i) - \bar{r}_i.
\end{equation}
The residual $\bar{r}_0$ is computed explicitly in Algorithm \ref{algo::pcg}, and the gap $f_0$ is thus the round-off from computing $\bar{r}_0$ from $A$, $\bar{x}_0$ and $b$, i.e., $f_0 = b-A\bar{x}_0 - \text{fl}(b-A\bar{x}_0)$. The norm of this initial gap is bounded by
\begin{equation}
	\|f_0\| \leq \left( (\mu\sqrt{n}+1) \, \|A\| \, \|\bar{x}_0\| + \|b\| \right) \epsilon.
\end{equation}
In iteration $i$ we obtain the following formula for the gap by substituting the recursions \eqref{eq:pxr}:
\begin{eqnarray}
	f_{i+1} &=& (b-A\bar{x}_{i+1}) - \bar{r}_{i+1} \notag \\
					&=& b-A(\bar{x}_i + \bar{\alpha}_i \bar{p}_i + \delta_i^x) - (\bar{r}_i - \bar{\alpha}_i A \bar{p}_i + \delta_i^r) \notag \\
					&=& f_i - A\delta_i^x - \delta_i^r. \label{eq:recur_f}
\end{eqnarray}
This recursive formulation relates the residual error $f_{i+1}$ in step $i$ to the previous residual error $f_i$.
By taking norms on both sides, we obtain an upper bound on $\|f_{i+1}\|$ in function of the previous gap $\|f_i\|$:
\begin{equation}
	\|f_{i+1}\| \leq \|f_i\| + \|A\delta_i^x + \delta_i^r\|
\end{equation}
where we can use the bounds \eqref{eq:delta_pxr} to further rewrite the right-hand side bound as
\begin{eqnarray} \label{eq:bounds_pcg1}
	 \|A\delta_i^x + \delta_i^r\| 
				&\leq& \|A\| \, \|\delta_i^x\| + \| \delta_i^r\| \notag \\
				&\leq& \left( \|A\| \, \|\bar{x}_i\| + 2 \, |\bar{\alpha}_i| \, \|A\| \, \|\bar{p}_i\| + \|\bar{r}_i\| + (\mu\sqrt{n}+2) \, |\bar{\alpha}_i| \, \|A\| \, \|\bar{p}_i\| \right)\epsilon \notag \\
				&=& 	 \left( \|A\| \, \|\bar{x}_i\| + (\mu\sqrt{n}+4) \, |\bar{\alpha}_i| \, \|A\| \, \|\bar{p}_i\| + \|\bar{r}_i\| \right) \epsilon \notag \\
				&:=&	 e_i^f \epsilon.
\end{eqnarray}
Hence, with the above definition of the upper bound factor $e_i^f$, we obtain 
\begin{equation}
	\|f_{i+1}\| \leq \|f_i\| + e_i^f \epsilon,
\end{equation}
which gives an upper bound on the norm of the gap between the true and recursive residual in any iteration of the method based on the gap in the previous iteration. 

The recurrence \eqref{eq:recur_f} implies that in the classical CG method the gap $f_{i+1}$ is the sum of local rounding errors, i.e.,
\begin{equation}\label{eq:gap_f0_cg}
	f_{i+1} = f_0  - \sum_{j=0}^i \left(A\delta_j^x + \delta_j^r\right).
\end{equation}
Hence, no amplification of local rounding errors occurs in classical CG, since \eqref{eq:gap_f0_cg} indicates that local rounding errors are simply accumulated, see also \cite{greenbaum1997estimating,gutknecht2000accuracy,sleijpen1994bicgstab}.

\subsection{Propagation of local rounding errors in Chronopoulos/Gear CG}

In so-called Chronopoulos/Gear CG (commonly denoted CG-CG in this manuscript), Algorithm \ref{algo::pcgcg}, an extra recurrence for the auxiliary variable $s_i$ is introduced, which in exact arithmetic equals $Ap_i$, and the auxiliary vectors $w_i = Au_i$ and $u_i = M^{-1}r_i$ are computed explicitly in each iteration, i.e., $\bar{w}_i = \text{fl}(A\bar{u}_i)$ and $\bar{u}_i = \text{fl}(M^{-1}\bar{r}_i)$. Alg.~\ref{algo::pcgcg} `avoids' communication by reducing the two global reduction phases of classical CG to one global synchronization (lines 11-12). The unpreconditioned version of the algorithm can be obtained by simply removing line 9 and replacing $u_i$ by $r_i$ where required. In finite precision arithmetic the corresponding recurrences in Alg.~\ref{algo::pcgcg} are
\begin{align}
	\bar{x}_{i+1} &= \bar{x}_i + \bar{\alpha}_i \bar{p}_i + \delta_i^x ,  &  
	\bar{p}_i 		&= \bar{u}_i + \bar{\beta}_i \bar{p}_{i-1} + \delta_i^p , \notag \\
	\bar{r}_{i+1} &= \bar{r}_i - \bar{\alpha}_i \bar{s}_i + \delta_i^r  , &
 	\bar{s}_i 		&= A\bar{u}_i + \bar{\beta}_i \bar{s}_{i-1} + \delta_i^s \label{eq:pxr_cgcg},
\end{align}
where the local rounding errors satisfy
\begin{eqnarray} 
	\|\delta_i^x\| &\leq& \left( \|\bar{x}_i\| + 2 \, |\bar{\alpha}_i| \, \|\bar{p}_i\| \right) \epsilon , \notag \\
	\|\delta_i^r\| &\leq& \left( \|\bar{r}_i\| + 2 \, |\bar{\alpha}_i| \, \|\bar{s}_i\| \right) \epsilon, \notag \\
	\|\delta_i^p\| &\leq& \left( (\tilde{\mu}\sqrt{n}+1) \, \|M^{-1}\| \, \|\bar{r}_i\| + 2 \, |\bar{\beta}_i| \, \|\bar{p}_{i-1}\| \right) \epsilon, \notag \\
	\|\delta_i^s\| &\leq& \left( (\mu\sqrt{n}+\tilde{\mu}\sqrt{n}+1) \, \|A\| \, \|M^{-1}\| \, \|\bar{r}_i\| + 2 \, |\bar{\beta}_i| \, \|\bar{s}_{i-1}\| \right) \epsilon \label{eq:delta_pxr_cgcg} ,
\end{eqnarray}
where $\tilde{\mu}$ is the maximum number of nonzeros in any row of the operator $M^{-1}$. To estimate the gap between the true and recursive residual we again substitute the recursions \eqref{eq:pxr_cgcg} in $f_i = (b - A\bar{x}_i) - \bar{r}_i$. Note that we have
\begin{equation} \label{eq:f0_cgcg}
	\|f_0\| \leq \left( (\mu\sqrt{n}+1) \, \|A\| \, \|\bar{x}_0\| + \|b\| \right) \epsilon,
\end{equation}
since the recursion for $\bar{x}_i$ in CG-CG is identical to that in CG, see \eqref{eq:pxr}. The gap in iteration $i$ is given by
\begin{eqnarray}
	f_{i+1} &=& (b-A\bar{x}_{i+1}) - \bar{r}_{i+1} \notag \\
					&=& b-A(\bar{x}_i + \bar{\alpha}_i \bar{p}_i + \delta_i^x) - (\bar{r}_i - \bar{\alpha}_i \bar{s}_i + \delta_i^r) \notag \\
					&=& f_i - \bar{\alpha}_i g_i - A\delta_i^x - \delta_i^r, \label{eq:recur_f_cgcg}
\end{eqnarray}
where $g_i = A \bar{p}_i - \bar{s}_i$, that is the gap between the explicit quantity $A\bar{p}_i$ and the recursively computed auxiliary variable $\bar{s}_i$. For the latter gap, it holds for $i = 0$ that
\begin{equation} \label{eq:g0_cgcg}
	\|g_0\| \leq \mu\sqrt{n} \, \|A\| \, \|\bar{p}_0\| \, \epsilon.
\end{equation}
Indeed, in iteration $i = 0$, the gap $g_0$ is the round-off on the explicit computation $A\bar{p}_0$, i.e., $g_0 = A\bar{p}_0 - \text{fl}(A\bar{p}_0)$. In iteration $i > 0$ the variable $\bar{s}_i$ is computed recursively, and it holds that
\begin{eqnarray}
	g_i &=& A \bar{p}_i - \bar{s}_i \notag \\
			&=& A (\bar{u}_i + \bar{\beta}_i \bar{p}_{i-1} +\delta_i^p) - (A\bar{u}_i + \bar{\beta}_i \bar{s}_{i-1} + \delta_i^s) \notag \\
			&=& \bar{\beta}_i g_{i-1} + A \delta_i^p - \delta_i^s. \label{eq:recur_g_cgcg}
\end{eqnarray}

\begin{algorithm}[t]
  \caption{Preconditioned Chronopoulos/Gear CG}
  \label{algo::pcgcg}
    \begin{algorithmic}[1]
    	\Procedure{prec-cg-cg}{$A$, $M^{-1}$, $b$, $x_0$}
      \State $r_0 := b - Ax_0$; $u_0 := M^{-1}r_0$; $w_0 := Au_0$
      \State $\alpha_0 := (r_0,u_0)/(w_0,u_0)$; $\beta_0:=0$; $\gamma_0 := (r_0,u_0)$
      \For{$i = 0, \dots$}
      \State $p_i := u_i + \beta_i p_{i-1}$
      \State {$s_i := w_i + \beta_i s_{i-1}$}
      \State {$x_{i+1} := x_i + \alpha_i p_i$}
      \State {$r_{i+1} := r_i - \alpha_i s_i$}
      \State {$u_{i+1} := M^{-1} r_{i+1}$}
      \State {$w_{i+1} := Au_{i+1}$}
      \State {$\gamma_{i+1} := (r_{i+1},u_{i+1})$}
      \State {$\delta := (w_{i+1},u_{i+1})$}
      \State {$\beta_{i+1} := \gamma_{i+1}/\gamma_i$}
      \State {$\alpha_{i+1} := (\delta/\gamma_{i+1}-\beta_{i+1}/\alpha_i)^{-1}$}
      \EndFor
			\EndProcedure
    \end{algorithmic}
\end{algorithm}

\noindent The residual gap is thus given by the coupled recursions
\begin{equation} \label{eq:pcgcg_system}
\begin{bmatrix}
 f_{i+1}  \\  g_i 
\end{bmatrix} = 
\begin{bmatrix}
    1 & -\bar{\alpha}_i \bar{\beta}_i \\ 0  & \bar{\beta}_i
\end{bmatrix}
\begin{bmatrix}
f_i \\ g_{i-1}
\end{bmatrix} +
\begin{bmatrix}
  - A\delta_i^x - \delta_i^r - \bar{\alpha}_i \left( A \delta_i^p - \delta_i^s \right) \\ A \delta_i^p - \delta_i^s
\end{bmatrix}.
\end{equation}
Taking norms, we obtain the upper bounds
\begin{equation} \label{eq:pcgcg_bounds}
\begin{bmatrix}
 \|f_{i+1}\|  \\  \|g_i\| 
\end{bmatrix} \leq
\begin{bmatrix}
    1 & |\bar{\alpha}_i \bar{\beta}_i| \\ 0  & |\bar{\beta}_i|
\end{bmatrix}
\begin{bmatrix}
\|f_i\| \\ \|g_{i-1}\|
\end{bmatrix} +
\begin{bmatrix}
  \| A\delta_i^x + \delta_i^r \| + | \bar{\alpha}_i | \| A \delta_i^p - \delta_i^s \| \\ \| A \delta_i^p - \delta_i^s \|
\end{bmatrix}.
\end{equation}
This bound can be further rewritten into more tractable expressions using the bounds for the local rounding errors in \eqref{eq:delta_pxr_cgcg}, i.e.,
\begin{eqnarray} \label{eq:pcgcg_bounds3}
	 \|A\delta_i^x + \delta_i^r\| 
				&\leq& \|A\| \, \|\delta_i^x\| + \| \delta_i^r\| \notag \\
				&\leq& \left( \|A\| \, \|\bar{x}_i\| + 2 \, |\bar{\alpha}_i| \, \|A\| \, \|\bar{p}_i\| + \|\bar{r}_i\| + 2 \, |\bar{\alpha}_i| \, \|\bar{s}_i\| \right)\epsilon \notag \\
				&:=&	 e_i^f \epsilon,
\end{eqnarray}
and
\begin{eqnarray}
	 \|A\delta_i^p - \delta_i^s\| 
				&\leq& \|A\| \, \|\delta_i^p\| + \| \delta_i^s\| \notag \\
				&\leq& \left( (\tilde{\mu}\sqrt{n}+1) \, \|A\| \, \|M^{-1}\| \, \|\bar{r}_i\| + 2 \, |\bar{\beta}_i| \, \|A\| \, \|\bar{p}_{i-1}\| \notag \right. \\ 
								&& ~ + \left. (\mu\sqrt{n}+\tilde{\mu}\sqrt{n}+1) \, \|A\| \, \|M^{-1}\| \, \|\bar{r}_i\| + 2 \, |\bar{\beta}_i| \, \|\bar{s}_{i-1}\| \right)\epsilon \notag \\
				&=& 	 \left( 2 \, |\bar{\beta}_i| \, \|A\| \, \|\bar{p}_{i-1}\| + ((\mu+2\tilde{\mu})\sqrt{n}+2) \, \|A\| \, \|M^{-1}\| \, \|\bar{r}_i\| \right. \notag \\
								&& ~ + \left. 2 \, |\bar{\beta}_i| \, \|\bar{s}_{i-1}\| \right) \epsilon \notag \\
				&:=&	 e_i^g \epsilon.
\end{eqnarray}
Note that the definitions of the bounds are local to each subsection; definition \eqref{eq:pcgcg_bounds3} above holds for CG-CG, and should not be confused with the earlier identical notation defined by \eqref{eq:bounds_pcg1} for the bound in classical CG.
Hence, with the factors $e_i^f$ and $e_i^g$ defined as above, the norm of the gap between the true and recursive residual is bounded by the recursively defined system of upper bounds
\begin{equation} \label{eq:pcgcg_bounds2}
\begin{bmatrix}
 \|f_{i+1}\|  \\  \|g_i\| 
\end{bmatrix} \leq
\begin{bmatrix}
    1 & |\bar{\alpha}_i \bar{\beta}_i| \\ 0  & |\bar{\beta}_i|
\end{bmatrix}
\begin{bmatrix}
\|f_i\| \\ \|g_{i-1}\|
\end{bmatrix} +
\begin{bmatrix}
  e_i^f \epsilon + | \bar{\alpha}_i | \, e_i^g \epsilon \\ e_i^g \epsilon
\end{bmatrix}.
\end{equation}
From \eqref{eq:pcgcg_system} it can be derived by induction that the residual gap in iteration $i$ is:
\begin{equation} \label{eq:gap_f0_cgcg} \vspace{-0.1cm}
	f_{i+1} = f_0 - \sum_{j=0}^i \left(A\delta_j^x + \delta_j^r\right) - \sum_{j=0}^i \bar{\alpha}_j \left[ \left(\prod_{k=1}^j \bar{\beta}_k\right) g_0 + \sum_{k=1}^j\left(\prod_{l=k+1}^j \bar{\beta}_l \right) \left(A\delta_k^p - \delta_k^s\right) \right].
\end{equation}
Note that this is in sharp contrast to the error behavior of the residual gap in the classical CG algorithm, where the gap after $i+1$ steps is a simple sum of local rounding errors, see \eqref{eq:gap_f0_cg}. Indeed, the local rounding errors $\left(A\delta_k^p - \delta_k^s\right)$ $(1 \leq k \leq j)$ that contribute to the difference $f_{i+1} = \left(b-A\bar{x}_{i+1}\right)-\bar{r}_{i+1}$ in \eqref{eq:gap_f0_cgcg} are potentially amplified by the factors $\prod_{l=k+1}^j \bar{\beta}_l$. Note that in exact arithmetic this product is
\begin{equation}
	\prod_{l=k+1}^j \beta_l  =  \frac{\| r_j \|^2}{\| r_k \|^2}, \quad 1 \leq k \leq j,
\end{equation}
which may be large for some $k \leq j$. Consequently, like the three-term recurrence CG algorithm \cite{stiefel1955relaxationsmethoden} which was analyzed in \cite{gutknecht2000accuracy}, see also \cite{carson2016numerical}, the CG-CG method may suffer from a dramatic amplification of local rounding errors throughout the algorithm, and the accuracy achieved by Alg.~\ref{algo::pcgcg} can be significantly worse compared to Alg.~\ref{algo::pcg}.

\subsection{Propagation of local rounding errors in pipelined CG}

In preconditioned pipelined CG, Algorithm \ref{algo::ppipe-cg}, additional recurrences are introduced for the auxiliary variables $w_i$, $z_i$, $u_i$ and $q_i$, which respectively equal $A u_i$, $A q_i$, $M^{-1} r_i$ and $M^{-1} s_i$ in exact arithmetic, whereas $v_i = A m_i$ and $m_i = M^{-1} w_i$ are computed explicitly. In addition to reducing communication, Alg.~\ref{algo::ppipe-cg} `hides' the communication phase (lines 4-5) behind the \textsc{spmv} and preconditioner application (lines 6-7). Replacing the recurrences by their finite precision equivalents, we have
\begin{align}
	\bar{x}_{i+1} &= \bar{x}_i + \bar{\alpha}_i \bar{p}_i + \delta_i^x, &
	\bar{p}_i 		&= \bar{u}_i + \bar{\beta}_i \bar{p}_{i-1} + \delta_i^p , \notag \\
	\bar{r}_{i+1} &= \bar{r}_i - \bar{\alpha}_i \bar{s}_i + \delta_i^r,&
	\bar{s}_i 		&= \bar{w}_i + \bar{\beta}_i \bar{s}_{i-1} + \delta_i^s, \notag \\
	\bar{w}_{i+1} &= \bar{w}_i - \bar{\alpha}_i \bar{z}_i + \delta_i^w, &
	\bar{z}_i     &= A\bar{m}_i + \bar{\beta}_i \bar{z}_{i-1} + \delta_i^z, \notag \\
	\bar{u}_{i+1} &= \bar{u}_i - \bar{\alpha}_i \bar{q}_i + \delta_i^u , &
	\bar{q}_i     &= \bar{m}_i + \bar{\beta}_i \bar{q}_{i-1} + \delta_i^q,\label{eq:pxr_pipecg}
\end{align}
where the local rounding errors are bounded by
\begin{eqnarray} 
	\|\delta_i^x\| &\leq& \left( \|\bar{x}_i\| + 2 \, |\bar{\alpha}_i| \, \|\bar{p}_i\| \right) \epsilon , \notag \\
	\|\delta_i^p\| &\leq& \left( \|\bar{u}_i\| + 2 \, |\bar{\beta}_i| \, \|\bar{p}_{i-1}\| \right) \epsilon , \notag \\
	\|\delta_i^r\| &\leq& \left( \|\bar{r}_i\| + 2 \, |\bar{\alpha}_i| \, \|\bar{s}_i\| \right) \epsilon , \notag \\
	\|\delta_i^s\| &\leq& \left( \|\bar{w}_i\| + 2 \, |\bar{\beta}_i| \, \|\bar{s}_{i-1}\| \right) \epsilon , \notag \\
	\|\delta_i^w\| &\leq& \left( \|\bar{w}_i\| + 2 \, |\bar{\alpha}_i| \, \|\bar{z}_i\| \right) \epsilon  , \notag \\
	\|\delta_i^z\| &\leq& \left( (\mu\sqrt{n}+\tilde{\mu}\sqrt{n}+1)  \, \|A\| \, \|M^{-1}\| \, \|\bar{w}_i\| + 2 \, |\bar{\beta}_i| \, \|\bar{z}_{i-1}\| \right) \epsilon  , \notag \\
	\|\delta_i^u\| &\leq& \left( \|\bar{u}_i\| + 2 \, |\bar{\alpha}_i| \, \|\bar{q}_i\| \right) \epsilon  , \notag \\
	\|\delta_i^q\| &\leq& \left( (\tilde{\mu}\sqrt{n}+1) \, \|M^{-1}\| \, \|\bar{w}_i\| + 2 \, |\bar{\beta}_i| \, \|\bar{q}_{i-1}\| \right) \epsilon   \label{eq:delta_pxr_pipecg} .
\end{eqnarray}

\begin{algorithm}[t]
  \caption{Preconditioned pipelined CG}
  \label{algo::ppipe-cg}
  \begin{algorithmic}[1]
  	\Procedure{prec-p-cg}{$A$, $M^{-1}$, $b$, $x_0$}
    \State $r_0 := b - Ax_0$; $u_0:= M^{-1} r_0$; $w_0 := Au_0$
    \For{$i = 0,\dots$}
    \State $\gamma_i :=(r_i,u_i)$
    \State $\delta := (w_i,u_i)$
    \State $m_i := M^{-1} w_i$
    \State $v_i := A m_i$
    \If{$i>0$}
    \State $\beta_i := \gamma_i/\gamma_{i-1}$; $\alpha_i := (\delta/\gamma_i - \beta_i/\alpha_{i-1})^{-1}$
    \Else
    \State $\beta_i := 0$; $\alpha_i := \gamma_i/\delta$
    \EndIf
    \State $z_i := v_i + \beta_i z_{i-1}$
    \State $q_i := m_i + \beta_i q_{i-1}$
    \State $s_i := w_i + \beta_i s_{i-1}$
    \State $p_i := u_i + \beta_i p_{i-1}$
    \State $x_{i+1} := x_i + \alpha_i p_i$
    \State $r_{i+1} := r_i - \alpha_i s_i$
    \State $u_{i+1} := u_i - \alpha_i q_i$
    \State $w_{i+1} := w_i - \alpha_i z_i$
    \EndFor
    \EndProcedure
  \end{algorithmic}
\end{algorithm}

The gap $f_i = (b-A\bar{x}_i) - \bar{r}_i$ can then be calculated similarly to \eqref{eq:recur_f_cgcg}. The initial gap $f_0$ satisfies \eqref{eq:f0_cgcg}, and in iteration $i$ we have
\begin{eqnarray}
	f_{i+1} &=& (b-A\bar{x}_{i+1}) - \bar{r}_{i+1} \notag \\
					&=& b-A(\bar{x}_i + \bar{\alpha}_i \bar{p}_i + \delta_i^x) - (\bar{r}_i - \bar{\alpha}_i \bar{s}_i + \delta_i^r) \notag \\
					&=& f_i - \bar{\alpha}_i g_i - A\delta_i^x - \delta_i^r. \label{eq:recur_f_pipecg}
\end{eqnarray}
The residual gap is again coupled to $g_i = A\bar{p}_i - \bar{s}_i$, which can be written as
\begin{eqnarray}
	g_i &=& A \bar{p}_i - \bar{s}_i \notag \\
			&=& A (\bar{u}_i + \bar{\beta}_i \bar{p}_{i-1} +\delta_i^p) - (\bar{w}_i + \bar{\beta}_i \bar{s}_{i-1} + \delta_i^s) \notag \\
			&=& h_i + \bar{\beta}_i g_{i-1} + A \delta_i^p - \delta_i^s, \label{eq:recur_g_pipecg}
\end{eqnarray}
where $g_0$ satisfies \eqref{eq:g0_cgcg} and we define $h_i = A\bar{u}_i - \bar{w}_i$. Instead of being computed explicitly, the auxiliary variable $\bar{w}_i$ is also computed recursively in pipelined CG, leading to an additional coupling of the residual gap $f_i$ to the difference $h_i$. For $i = 0$, it holds that the norm of the gap $h_i$ is bounded by
\begin{equation} \label{eq:h0_pipecg}
	\|h_0\| \leq \mu\sqrt{n} \, \|A\| \, \|\bar{u}_0\| \, \epsilon.
\end{equation}
Substituting the recurrences \eqref{eq:pxr_pipecg}, we find that the gap between $A \bar{u}_{i+1}$ and $\bar{w}_{i+1}$ in iteration $i$ is
\begin{eqnarray}
	h_{i+1} &=& A \bar{u}_{i+1} - \bar{w}_{i+1} \notag \\
					&=& A (\bar{u}_i - \bar{\alpha}_i \bar{q}_i +\delta_i^u) - (\bar{w}_i - \bar{\alpha}_i \bar{z}_i + \delta_i^w) \notag \\
					&=& h_i - \bar{\alpha}_i j_i + A \delta_i^u - \delta_i^w, \label{eq:recur_h_pipecg}
\end{eqnarray}
which relates the residual error to the error $j_i = A\bar{q}_i - \bar{z}_i$. The latter gap is due to the recursive computation of the auxiliary variable $\bar{z}_i$. For $i=0$, we can bound the norm of $j_i$ by the norm of the round-off, i.e.,
\begin{equation} \label{eq:j0_pipecg}
	\|j_0\| \leq \mu\sqrt{n} \, \|A\| \, \|\bar{q}_0\| \, \epsilon.
\end{equation}
Using again the recursive definitions \eqref{eq:pxr_pipecg}, we obtain
\begin{eqnarray}
	j_i &=& A \bar{q}_i - \bar{z}_i \notag \\
			&=& A (\bar{m}_i + \bar{\beta}_i \bar{q}_{i-1} + \delta_i^q) - (A\bar{m}_i + \bar{\beta}_i \bar{z}_{i-1} + \delta_i^z) \notag \\
			&=& \bar{\beta}_i j_{i-1} + A \delta_i^q - \delta_i^z. \label{eq:recur_j_pipecg}
\end{eqnarray}
Hence, for pipelined CG, the residual gap is given by the system of coupled equations
\begin{equation} \label{eq:pipecg_system}
\begin{bmatrix}
 f_{i+1}  \\  g_i \\ h_{i+1} \\ j_i
\end{bmatrix} = 
\begin{bmatrix}
    1 & -\bar{\alpha}_i \bar{\beta}_i & -\bar{\alpha}_i & 0 \\ 0 & \bar{\beta}_i & 1 & 0 \\ 0 & 0 & 1 & -\bar{\alpha}_i \bar{\beta}_i \\ 0 & 0 & 0 & \bar{\beta}_i
\end{bmatrix}
\begin{bmatrix}
f_i \\ g_{i-1} \\ h_i \\ j_{i-1}
\end{bmatrix} +
\begin{bmatrix}
  - A\delta_i^x - \delta_i^r - \bar{\alpha}_i \left( A \delta_i^p - \delta_i^s \right) \\ 
	A \delta_i^p - \delta_i^s \\
	A\delta_i^u - \delta_i^w - \bar{\alpha}_i \left( A \delta_i^q - \delta_i^z \right) \\ 
	A \delta_i^q - \delta_i^z 
\end{bmatrix}.
\end{equation}
Taking norms of both sides in \eqref{eq:pipecg_system}, we arrive at the following coupled system of upper bounds for the gaps in pipelined CG:
\begin{equation} \label{eq:pipecg_bounds}
\begin{bmatrix}
 \|f_{i+1} \| \\  \|g_i\| \\ \|h_{i+1}\| \\ \|j_i\|
\end{bmatrix} \leq 
\begin{bmatrix}
    1 & |\bar{\alpha}_i \bar{\beta}_i| & |\bar{\alpha}_i| & 0\\  0& |\bar{\beta}_i| & 1 &0 \\ 0& 0& 1 & |\bar{\alpha}_i \bar{\beta}_i| \\ 0& 0& 0& |\bar{\beta}_i|
\end{bmatrix}
\begin{bmatrix}
	\|f_i\| \\ \|g_{i-1}\| \\ \|h_i\| \\ \|j_{i-1}\|
\end{bmatrix} +
\begin{bmatrix}
  \| A\delta_i^x + \delta_i^r \| + |\bar{\alpha}_i| \| A \delta_i^p - \delta_i^s \| \\ 
	\| A \delta_i^p - \delta_i^s \| \\
	\| A\delta_i^u - \delta_i^w \| + |\bar{\alpha}_i| \| A \delta_i^q - \delta_i^z \| \\ 
	\| A \delta_i^q - \delta_i^z \|
\end{bmatrix}.
\end{equation}
The per-iteration additions on the right-hand side in \eqref{eq:pipecg_bounds} can be bounded further using the local error bounds \eqref{eq:delta_pxr_pipecg}. 
We hence obtain the computable bounds
\begin{eqnarray}
	 \|A\delta_i^x + \delta_i^r\| 
				&\leq& \|A\| \, \|\delta_i^x\| + \| \delta_i^r\| \notag \\
				&\leq& \left( \|A\| \, \|\bar{x}_i\| + 2 \, |\bar{\alpha}_i| \, \|A\| \, \|\bar{p}_i\| + \|\bar{r}_i\| + 2 \, |\bar{\alpha}_i| \, \|\bar{s}_i\| \right)\epsilon \notag \\
				&:=&	 e_i^f \epsilon , \label{eq:eif}
\end{eqnarray}
\begin{eqnarray}
	 \|A\delta_i^p - \delta_i^s\| 
				&\leq& \|A\| \, \|\delta_i^p\| + \| \delta_i^s\| \notag \\
				&\leq& \left( \|A\| \, \|\bar{u}_i\| + 2 \, |\bar{\beta}_i| \, \|A\| \, \|\bar{p}_{i-1}\| + \|\bar{w}_i\| + 2 \, |\bar{\beta}_i| \, \|\bar{s}_{i-1}\| \right) \epsilon \notag \\
				&:=&	 e_i^g \epsilon , \label{eq:eig}
\end{eqnarray}
\begin{eqnarray}
	 \|A\delta_i^u - \delta_i^w\| 
				&\leq& \|A\| \, \|\delta_i^u\| + \| \delta_i^w\| \notag \\
				&\leq& \left( \|A\| \, \|\bar{u}_i\| + 2 \, |\bar{\alpha}_i| \, \|A\| \, \|\bar{q}_i\| + \|\bar{w}_i\| + 2 \, |\bar{\alpha}_i| \, \|\bar{z}_i\| \right)\epsilon \notag \\
				&:=&	 e_i^h \epsilon , \label{eq:eih}
\end{eqnarray}
\begin{eqnarray}
	 \|A\delta_i^q - \delta_i^z\| 
				&\leq& \|A\| \, \|\delta_i^q\| + \| \delta_i^z\| \notag \\
				&\leq& \left( (\tilde{\mu}\sqrt{n}+1) \, \|A\| \, \|M^{-1}\| \, \|\bar{w}_i\| + 2 \, |\bar{\beta}_i| \, \|A\| \, \|\bar{q}_{i-1}\| \notag \right. \\ 
								&& ~ + \left. (\mu\sqrt{n}+\tilde{\mu}\sqrt{n}+1) \, \|A\| \, \|M^{-1}\| \, \|\bar{w}_i\| + 2 \, |\bar{\beta}_i| \, \|\bar{z}_{i-1}\| \right)\epsilon \notag \\
				&=& 	 \left( 2 \, |\bar{\beta}_i| \, \|A\| \, \|\bar{q}_{i-1}\| + ((\mu+2\tilde{\mu})\sqrt{n}+2) \, \|A\| \, \|M^{-1}\| \, \|\bar{w}_i\| \right. \notag \\
								&& ~ + \left. 2 \, |\bar{\beta}_i| \, \|\bar{z}_{i-1}\| \right) \epsilon \notag \\
				&:=&	 e_i^j \epsilon . \label{eq:eij}
\end{eqnarray}
This yields the following system of upper bounds on the norm of the gap between the true and recursive residual:
\begin{equation} \label{eq:pipecg_bounds2}
\begin{bmatrix}
 \|f_{i+1} \| \\  \|g_i\| \\ \|h_{i+1}\| \\ \|j_i\|
\end{bmatrix} \leq 
\begin{bmatrix}
    1 & |\bar{\alpha}_i \bar{\beta}_i| & |\bar{\alpha}_i| & 0 \\ 0 & |\bar{\beta}_i| & 1 & 0\\ 0& 0& 1 & |\bar{\alpha}_i \bar{\beta}_i| \\ 0& 0& 0& |\bar{\beta}i|
\end{bmatrix}
\begin{bmatrix}
	\|f_i\| \\ \|g_{i-1}\| \\ \|h_i\| \\ \|j_{i-1}\|
\end{bmatrix} +
\begin{bmatrix}
  e_i^f \epsilon + |\bar{\alpha}_i| e_i^g \epsilon \\ 
	e_i^g \epsilon \\
	e_i^h \epsilon + |\bar{\alpha}_i| e_i^j \epsilon \\ 
	e_i^j \epsilon
\end{bmatrix},
\end{equation}
where $e_i^f$, $e_i^g$, $e_i^h$ and $e_i^j$ are defined in \eqref{eq:eif}-\eqref{eq:eij}. 

By induction, \eqref{eq:pipecg_system} can be reformulated into an expression for the residual gap $f_{i+1}$ with respect to $f_0, g_0, h_0,j_0$ and local rounding errors, similar to \eqref{eq:gap_f0_cgcg}, i.e.,
\begin{equation} \label{eq:gap_f0_pcg}
	f_{i+1} = f_0 - \sum_{j=0}^i \bar{\alpha}_j g_j - \sum_{j=0}^i \left(A\delta_j^x + \delta_j^r\right),
\end{equation}
where
\begin{equation} \label{eq:gap_g0_pcg}
	g_j = \left(\prod_{k=1}^j\bar{\beta}_k\right) g_0 + \sum_{k=1}^j \left(\prod_{l=k+1}^j \bar{\beta}_l\right) \left(A\delta_k^p-\delta_k^s\right) + \sum_{k=1}^j \left(\prod_{l=k+1}^j \bar{\beta}_l\right) h_k,
\end{equation}
with
\begin{equation} \label{eq:gap_h0_pcg}
h_k = h_0 - \sum_{l=0}^{k-1} \bar{\alpha}_l j_l + \sum_{l=0}^{k-1} \left(A \delta_l^u - \delta_l^w\right),
\end{equation}
and where
\begin{equation} \label{eq:gap_j0_pcg}
j_l = \left( \prod_{m=1}^l \bar{\beta}_m\right)	j_0 + \sum_{m=1}^l \left(\prod_{n=m+1}^l \bar{\beta}_n\right) \left(A\delta_m^q - \delta_m^z\right).
\end{equation}
It is clear from \eqref{eq:gap_f0_pcg}-\eqref{eq:gap_j0_pcg} that, due to the extra recurrence relations, the propagation of local rounding errors may be even more dramatic for p-CG, since $\left(A\delta_k^p-\delta_k^s\right)$ $(1 \leq k \leq j)$, $\left(A \delta_l^u - \delta_l^w\right)$ $(0\leq l \leq k-1)$ and $\left(A\delta_m^q - \delta_m^z\right)$ $(1\leq m\leq l)$ are all potentially amplified during the algorithm. This may lead to significantly reduced maximal attainable accuracy compared to both classical CG and CG-CG.

Note that the auxiliary variables $u_i$ and $q_i$, which in exact arithmetic represent the preconditioned versions of the residual $r_i$ and the auxiliary variable $s_i$ respectively, are also computed recursively in pipelined CG. Hence, we can write down an analogous derivation for the gap between the explicit and recursive preconditioned residual, that is, $k_i = M^{-1} \bar{r}_i - \bar{u}_i$. For $i = 0$ we have
\begin{equation} \label{eq:k0_pipecg}
	\|k_0\| \leq \tilde{\mu}\sqrt{n} \, \|M^{-1}\| \, \|\bar{r}_0\| \, \epsilon,
\end{equation}
and in iteration $i$ we find
\begin{eqnarray}
	k_{i+1} &=& M^{-1} \bar{r}_{i+1} - \bar{u}_{i+1} \notag \\
					&=& M^{-1} (\bar{r}_i - \bar{\alpha}_i \bar{s}_i + \delta_i^r) - (\bar{u}_i - \bar{\alpha}_i \bar{q}_i + \delta_i^u) \notag \\
					&=& k_i - \bar{\alpha}_i \ell_i + M^{-1}\delta_i^r - \delta_i^u, \label{eq:recur_k_pipecg}
\end{eqnarray}
where we define $\ell_i = M^{-1} \bar{s}_i - \bar{q}_i$. Finally, for the gap between the explicit quantity $M^{-1} \bar{s}_i$ and the recursive variable $\bar{q}_i$, we have for $i = 0$ that
\begin{equation} \label{eq:l0_pipecg}
	\|\ell_0\| \leq \tilde{\mu}\sqrt{n} \, \|M^{-1}\| \, \|\bar{s}_0\| \, \epsilon.
\end{equation}
By once again inserting the recurrences from \eqref{eq:pxr_pipecg}, we find that $\ell_i$ in iteration $i$ is
\begin{eqnarray}
	\ell_i &=& M^{-1} \bar{s}_i - \bar{q}_i \notag \\
			&=& M^{-1} (\bar{w}_i + \bar{\beta}_i \bar{s}_{i-1} + \delta_i^s) - (\bar{m}_i + \bar{\beta}_i \bar{q}_{i-1} + \delta_i^q) \notag \\
			&=& \bar{\beta}_i \ell_{i-1} + M^{-1}\delta_i^s - \delta_i^q. \label{eq:recur_l_pipecg}
\end{eqnarray}
The last equation in \eqref{eq:recur_l_pipecg} holds since $\bar{m}_i$ is computed explicitly as $M^{-1} \bar{w}_i$ in Algorithm \ref{algo::ppipe-cg}.
This leads to a separate system of coupled recurrences for the gap on the preconditioned residual $k_{i+1}$:
\begin{equation} \label{eq:pipecg_system2}
\begin{bmatrix}
 k_{i+1}  \\  \ell_i 
\end{bmatrix} = 
\begin{bmatrix}
    1 & -\bar{\alpha}_i \bar{\beta}_i \\ 0  & \bar{\beta}_i
\end{bmatrix}
\begin{bmatrix}
	k_i \\ \ell_{i-1}
\end{bmatrix} +
\begin{bmatrix}
  M^{-1}\delta_i^r - \delta_i^u - \bar{\alpha}_i \left( M^{-1} \delta_i^s - \delta_i^q \right) \\ M^{-1} \delta_i^s - \delta_i^q
\end{bmatrix}.
\end{equation}
However, since the gap $k_{i+1}$ is uncoupled from the residual gap $f_{i+1}$, the coupled recurrences \eqref{eq:recur_k_pipecg}-\eqref{eq:recur_l_pipecg} are not be taken into account when establishing bounds for the norm of the residual gap in \eqref{eq:pipecg_bounds2}. 

\subsection{A practical estimate for the residual gap} \label{sec:practical}


In the previous sections we have derived upper bounds for the norm of the residual gap for several variants of the CG algorithm. Although insightful from an analytical perspective, these bounds are typically not sharp. Indeed, the bounds on the norms of the local rounding errors \eqref{eq:delta_pxr_pipecg} may largely overestimate the actual error norms, see the discussion in \cite{greenbaum1997estimating}, p.\,541. For example, the right-hand side of \eqref{eq:pipecg_bounds2} could be much larger than the left-hand side, and hence could provide a poor estimate for the actual residual gap. 

To obtain a more realistic estimate for the residual gap, we turn to statistical analysis of the rounding errors, see \cite{wilkinson1994rounding}. 
A well-known rule of thumb \cite{higham2002accuracy} is that a realistic error estimate can be obtained by replacing the dimension-dependent constants in a rounding error bound by their square root; thus if the bound is $f(n) \epsilon$, the error is approximately of order $\sqrt{f(n)}\epsilon$. Hence, instead of using the upper bounds \eqref{eq:eif}-\eqref{eq:eij}, we use the following approximations for the local error norms:
\begin{align}
\|A\delta_i^x + \delta_i^r\| &\approx \sqrt e_i^f \epsilon , \hspace{-1cm} & \hspace{-1cm}
\|A\delta_i^p - \delta_i^s\| &\approx \sqrt e_i^g \epsilon , \notag \\
\|A\delta_i^u + \delta_i^w\| &\approx \sqrt e_i^h \epsilon , \hspace{-1cm} & \hspace{-1cm}
\|A\delta_i^q - \delta_i^z\| &\approx \sqrt e_i^j \epsilon. 
\end{align}
Note that for the norms of the initial gaps in \eqref{eq:f0_cgcg}, \eqref{eq:g0_cgcg}, \eqref{eq:h0_pipecg} and \eqref{eq:j0_pipecg}, a similar square root rescaling of the respective dimension-dependent factors has to be applied. 
We hence assume that the norm of the residual gap in iteration $i$ of the pipelined CG algorithm can be estimated as follows:
\begin{equation} \label{eq:pipecg_estimate}
\begin{bmatrix}
 \|f_{i+1} \| \\  \|g_i\| \\ \|h_{i+1}\| \\ \|j_i\|
\end{bmatrix} \approx
\begin{bmatrix}
    1 & |\bar{\alpha}_i \bar{\beta}_i| & |\bar{\alpha}_i| & 0 \\ 0 & |\bar{\beta}_i| & 1 & 0\\ 0& 0& 1 & |\bar{\alpha}_i \bar{\beta}_i| \\ 0& 0& 0& |\bar{\beta}_i|
\end{bmatrix}
\begin{bmatrix}
	\|f_i\| \\ \|g_{i-1}\| \\ \|h_i\| \\ \|j_{i-1}\|
\end{bmatrix} +
\begin{bmatrix}
  \sqrt e_i^f \epsilon + |\bar{\alpha}_i| \sqrt e_i^g \epsilon \\ 
	\sqrt e_i^g \epsilon \\
	\sqrt e_i^h \epsilon + |\bar{\alpha}_i| \sqrt e_i^j \epsilon \\ 
	\sqrt e_i^j \epsilon
\end{bmatrix}.
\end{equation}
This approximation tends to yield a good (a posteriori) estimate for the actual residual gap, as illustrated by the numerical experiments in the next section.

A second practical remark concerns the computation of the matrix and preconditioner norms, $\|A\|$ and $\|M^{-1}\|$, in the estimate for the gap between the true and recursive residual. The use of the matrix 2-norm is often prohibited in practice, since it is computationally expensive for large scale systems. However, using the norm inequality $\| A \|_2 \leq \sqrt{n} \, \| A \|_\infty$, the matrix 2-norms in the estimate can be replaced by their respective maximum norms multiplied by $\sqrt{n}$. This slightly worsens the estimate, but provides practically computable quantities for $e_i^f$, $e_i^g$, $e_i^h$ and $e_i^j$. Alternatively, in the context of matrix-free computations, randomized probabilistic techniques for matrix norm computation may be used, see for example \cite{halko2011finding}.

A related issue concerns the computation of the norm of the preconditioner. The operator $M^{-1}$ is often not available in matrix form. This is the case when preconditioning the system with e.g., an Incomplete Cholesky factorization (ICC) type scheme, or any (stencil-based) scheme where $M^{-1}$ is not explicitly formed. For these commonly used preconditioning methods the norm $\| M^{-1} \|$ is unavailable. Explicit use of the preconditioner norm can be avoided by reformulating the local rounding error bounds \eqref{eq:delta_pxr_pipecg} with respect to the preconditioned variable $\bar{m}_i = M^{-1} \bar{w}_i$, that is:
\begin{eqnarray} 
	\|\delta_i^z\| &\leq& \left( (\mu\sqrt{n}+1)  \, \|A\| \, \|\bar{m}_i\| + 2 \, |\bar{\beta}_i| \, \|\bar{z}_{i-1}\| \right) \epsilon  , \notag \\
	\|\delta_i^q\| &\leq& \left( \|\bar{m}_i\| + 2 \, |\bar{\beta}_i| \, \|\bar{q}_{i-1}\| \right) \epsilon   \label{eq:delta_pxr_pipecg2} .
\end{eqnarray}
These bounds do not explicitly take the rounding error of the multiplication $M^{-1} \bar{w}_i$ into account, but rather implicitly bound the local rounding error using $\|\bar{m}_i\|$. With these local rounding error bounds, $e_i^j$ can now be defined analogous to \eqref{eq:eij} as
\begin{eqnarray}
	 \|A\delta_i^q - \delta_i^z\| 
				&\leq& \|A\| \, \|\delta_i^q\| + \| \delta_i^z\| \notag \\
				&\leq&  \left( \|A\| \, \|\bar{m}_i\| + 2 \, |\bar{\beta}_i| \, \|A\| \, \|\bar{q}_{i-1}\| \right. \notag \\
					&& ~ + \left. (\mu\sqrt{n}+1) \, \|A\| \, \|\bar{m}_i\| \, + 2 \, |\bar{\beta}_i| \, \|\bar{z}_{i-1}\| \right) \epsilon \notag \\
				&=& \left( (\mu\sqrt{n}+2) \, \|A\| \, \|\bar{m}_i\| + 2 \, |\bar{\beta}_i| \, \|A\| \, \|\bar{q}_{i-1}\| + 2 \, |\bar{\beta}_i| \, \|\bar{z}_{i-1}\| \right)\epsilon \notag \\
				&:=&	 e_i^j \epsilon . \label{eq:eij2}
\end{eqnarray}
With $e_i^j$ defined as in \eqref{eq:eij2}, \eqref{eq:pipecg_estimate} can also be used to estimate the residual gap in pipelined CG when the preconditioning matrix $M^{-1}$ is not formed explicitly. 

We point out that, as an alternative to computing a residual gap estimate, 
one could explicitly compute the residual $b-A\bar{x}_i$ in each iteration of the algorithm to keep track of the residual gap. However, calculating $b-A\bar{x}_i$ in each iteration is computationally much too expensive in practice. The computation of the estimate requires only the computation of $e_i^f$, $e_i^g$, $e_i^h$ and $e_i^j$ in each iteration of the algorithm, see \eqref{eq:pipecg_estimate}. Note that the scalar expressions \eqref{eq:eif}-\eqref{eq:eij} and \eqref{eq:eij2} contain several norms that are by default not computed in Algorithm \ref{algo::ppipe-cg}. However, these norm computations cause no additional communication overhead since they can be combined with the existing global communication phase on lines 4-5. Hence, the computational cost of calculating the residual gap estimate is negligible compared to computing the residual explicitly, and, in contrast to the true residual, the estimate \eqref{eq:pipecg_estimate} can be computed in real time in each iteration of the algorithm without additional computational overhead.

Finally, we note that the norms of some auxiliary variables, notably $\bar{p}_i$, $\bar{s}_i$, $\bar{q}_i$, $\bar{z}_i$ and $\bar{m}_i$, are unavailable in step $i$, since these vectors are not defined until \emph{after} the global reduction phase. Hence, their norms can be computed at the earliest in the global reduction phase of iteration $i+1$.
Consequently, in practical implementations the estimated norm $\|f_{i+1}\|$, defined by \eqref{eq:pipecg_estimate}, can only be computed in iteration $i+1$. This means that when including the estimates for the residual gap into the pipelined CG algorithm, a delay of one iteration on the estimates is unavoidable.

\section{Pipelined CG with automated residual replacement} \label{sec:pipelined}
	
In this section we propose an automated residual replacement strategy for pipelined CG, 
based on the estimate for the gap between the true and recursive residual proposed in Section \ref{sec:practical}. 
Although the derivation of our replacement strategy is partially based on heuristics and would certainly benefit from further theoretical investigation, 
this ad hoc countermeasure aims to improve the possibly dramatically reduced maximal attainable accuracy of the pipelined method.
The idea of performing manual residual replacements to increase the attainable accuracy of pipelined CG 
was already suggested in the original paper \cite{ghysels2014hiding}.
However, to establish an automated replacement strategy, a practically computable criterion for replacement should be available.
We suggest such a criterion under the assumption that the 
orthogonality of the Krylov basis vectors is not (critically) affected by rounding errors.
Although this condition may not always hold in practice, it enables to design a heuristic with low computational overhead that proves effective in many situations as described in Section \ref{sec:numerical}.

We follow the basic idea of residual replacement in Krylov subspace methods as discussed by 
Van der Vorst et al.~\cite{van2000residual} and Sleijpen et al.~\cite{sleijpen1996reliable,sleijpen2001differences}. 
In specific iterations of the algorithm, 
the vectors $\bar{r}_{i+1}$, $\bar{w}_{i+1}$, $\bar{u}_{i+1}$, $\bar{s}_i$, $\bar{z}_i$ and $\bar{q}_i$
which are computed recursively in iteration $i$, are instead computed explicitly, such that
\begin{align}
	\bar{r}_{i+1} &= \text{fl}(b-A \bar{x}_{i+1}), & \bar{u}_{i+1} &= \text{fl}(M^{-1} \bar{r}_{i+1}), & \bar{w}_{i+1} &= \text{fl}(A \bar{u}_{i+1}), \notag \\
	\bar{s}_i &= \text{fl}(A \bar{p}_i), & \bar{q}_i &= \text{fl}(M^{-1} \bar{s}_i) & \bar{z}_i &= \text{fl}(A \bar{q}_i). \label{eq:rrs}
\end{align}
Note how the current iterate $\bar{x}_{i+1}$ and the search direction $\bar{p}_i$ are evidently not replaced, 
since no explicit formulae for these vectors are available.

Two important caveats arise when incorporating a residual replacement in an iterative method. 
First, one could inquire if such a drastic replacement strategy does not destroy (or delay) convergence. 
A second, related question concerns the use of a criterion for the iteration in which replacements should take place. Since each residual 
replacement step comes at an additional cost of computing the \textsc{spmv}s in \eqref{eq:rrs}, an accurate criterion to determine 
the need for residual replacement that does not overestimate the total number of 
replacements is essential. 

We briefly recapitulate the main results from \cite{tong2000analysis} and \cite{van2000residual} below.
The recurrences for $\bar{r}_{i+1}$ and $\bar{p}_{i+1}$ in the (unpreconditioned) Algorithm \ref{algo::pcg} are 
\begin{align} \label{eq:recur_r_p}
  \bar{r}_{i+1} &= \bar{r}_i - \bar{\alpha}_i A\bar{p}_i + \delta^r_i , \notag \\
  \bar{p}_{i+1} &= \bar{r}_{i+1} + \bar{\beta}_{i+1} \bar{p}_{i} + \delta^p_i ,
\end{align}
where $\delta_i^r$ and $\delta_i^p$ are bounded by \eqref{eq:delta_pxr}. Combining the above recursions yields the perturbed Lanczos relation
\begin{equation} \label{eq:arnoldi}
	A Z_i = Z_i T_i - \frac{\|\bar{r}_{0} - \bar{\alpha}_{0} A\bar{p}_{0} \|}{\bar{\alpha}_i\|\bar{r}_1\| \|\bar{r}_i\|} \, \bar{r}_{i+1} e_i^T + F_i \quad \text{with} \quad Z_i = \left[\frac{\bar{r}_1}{\|\bar{r}_1\|},\ldots,\frac{\bar{r}_i}{\|\bar{r}_i\|}\right],
\end{equation}
see \cite{van2000residual}, where $T_i$ is a tridiagonal matrix and $F_i$ is a perturbation caused by the local rounding errors.
We assume that $Z_i$ is full rank (which might not be true in practice, see below).
A key result from \cite{tong2000analysis} states that if $\bar{r}_i$ satisfies the relation \eqref{eq:arnoldi}, then
\begin{equation} \label{eq:arnoldi_bound}
	\|\bar{r}_{i+1}\| \leq C_i \min_{p \in \mathcal{P}_i, p(0)=1} \| p(A-F_i Z^+_i) \, \bar{r}_1 \|,
\end{equation}
where $C_i > 0$ is an iteration-dependent constant. This result suggests that even if the 
perturbation $F_i$ is significantly larger than $\epsilon$, which is the case after residual replacement, 
the norm of the residual may not be significantly affected, as illustrated by the experiments presented in \cite{tong2000analysis}. Based on the bound 
\eqref{eq:arnoldi_bound}, Van der Vorst et al.~propose in \cite{van2000residual} to explicitly update the residuals and other vectors only when 
the residual norm is sufficiently large compared to the norm of the residual gap.
Performing replacements when $\|\bar{r}_i\|$ is small is not recommended, since this could negatively affect convergence. 
Note that although this exposition provides a useful intuition on when to perform residual replacements, it can in general not
be considered as a full theoretical validation of the latter, since the assumption that $Z_i$ is full rank is often not satisfied in practical computations.

A replacement in step $i$ eliminates the residual gap $f_{i+1}$. However, it should be carried out before $\|f_i\|$ becomes 
too large relative to $\|\bar{r}_i\|$. In analogy to \cite{van2000residual}, we use a threshold $\tau$, 
typically chosen as $\tau = \sqrt{\epsilon}$, and perform a residual replacement in step $i$ of Algorithm \ref{algo::ppipe-cg} if
\begin{equation}\label{eq:criterion}
	\| f_{i-1} \| \leq \tau \|\bar{r}_{i-1}\|  \quad  \text{and}  \quad  \| f_{i} \| > \tau \|\bar{r}_{i}\|.
\end{equation}
This criterion ensures that replacements are only allowed when $\|\bar{r}_i\|$ is sufficiently large with respect to $\|f_i\|$. 
Furthermore, no excess replacement steps are performed, thus keeping the total computational cost as low as possible. The residual gap model 
\eqref{eq:pipecg_estimate} allows for the practical implementation of the replacement criterion \eqref{eq:criterion} 
in pipelined CG, see Algorithm \ref{algo::prec-p-cg-rr}. Note that criterion \eqref{eq:criterion} does not compare the estimate 
$\| f_{i+1} \|$ to the current residual norm $\|\bar{r}_{i+1}\|$ computed in iteration $i$, since the norms of both 
these quantities are not yet available, see the discussion near the end of Section \ref{sec:practical}. This implies 
that some additional storage is required for these auxiliary variables between subsequent iterations. The resulting algorithm is called pipelined 
CG with automated residual replacement (p-CG-rr) and is detailed in Algorithm \ref{algo::prec-p-cg-rr}.

In practical computations conditioning of the Lanczos vectors matrix $Z_i$ in \eqref{eq:arnoldi} may be poor due to numerical loss of orthogonality, 
which may cause the pseudo-inverse $Z^+_i$ in \eqref{eq:arnoldi_bound} to become very ill-conditioned. This somewhat restricts the practical validity of the above argument. 
We again stress that the residual replacement strategy proposed in this section should be interpreted primarily as a practically useable heuristic that allows to improve accuracy, rather than a general and theoretically justified countermeasure to the loss of attainable accuracy.
Moreover, we point out that in finite precision arithmetic the residual replacement strategy itself may be a source of rounding errors in the algorithm \cite{carson2016numerical},
possibly leading to delayed convergence, see \cite{strakovs2002error}. An in-depth theoretical analysis of the effect of replacements on convergence is however beyond the scope of this work.

\begin{algorithm}[H]
  \caption{Preconditioned pipelined CG with automated residual replacement}
  \label{algo::prec-p-cg-rr}
  \begin{algorithmic}[1]
		{\small
  	\Procedure{prec-p-cg-rr}{$A$, $M^{-1}$, $b$, $x_0$}
    \State $r_0 := b - Ax_0$; $u_0:= M^{-1} r_0$; $w_0 := Au_0$; $\zeta := \|b\|_2$; $\tau := \sqrt{\epsilon}$
		\State $n = \text{length}(b)$; $\theta = \sqrt{n}\|A\|_{\infty}$; $\mu = \max(\text{rowsums}(A))$; replace := false
    \For{$i = 0,\dots$}
    \State $\gamma_i := (r_i,u_i)$; $\delta := (w_i,u_i)$; $\rho_{i+1} := \|r_i\|_2$
		\If{$i>0$}
		\State $\chi_i := \|x_{i-1}\|_2$; $\pi_i := \|p_{i-1}\|_2$; $\sigma_i := \|s_{i-1}\|_2$; $\xi_i := \|u_{i-1}\|_2$
		\State $\omega_i := \|w_{i-1}\|_2$; $\phi_i := \|q_{i-1}\|_2$; $\psi_i := \|z_{i-1}\|_2$; $\nu_i := \|m_{i-1}\|_2$
		\EndIf
    \State $m_i := M^{-1} w_i$
    \State $v_i := A m_i$
    \If{$i>0$}
    \State $\beta_i := \gamma_i/\gamma_{i-1}$; $\alpha_i := (\delta/\gamma_i - \beta_i/\alpha_{i-1})^{-1}$
    \Else
    \State $\beta_i := 0$; $\alpha_i := \gamma_i/\delta$
    \EndIf
    \State $z_i := v_i + \beta_i z_{i-1}$
    \State $q_i := m_i + \beta_i q_{i-1}$
    \State $s_i := w_i + \beta_i s_{i-1}$
    \State $p_i := u_i + \beta_i p_{i-1}$
    \State $x_{i+1} := x_i + \alpha_i p_i$
    \State $r_{i+1} := r_i - \alpha_i s_i$
    \State $u_{i+1} := u_i - \alpha_i q_i$
    \State $w_{i+1} := w_i - \alpha_i z_i$
		\If{$i>0$}
		\State $e_{i-1}^f := \theta \chi_i + 2 \left|\alpha_{i-1}\right| \theta \pi_i + \rho_i + 2 \left|\alpha_{i-1}\right| \sigma_i$
		\State $e_{i-1}^h := \theta \xi_i + 2 \left|\alpha_{i-1}\right| \theta \phi_i + \omega_i + 2 \left|\alpha_{i-1}\right| \psi_i$
		\If{$i>1$} 
		\State $e_{i-1}^g := \theta \xi_i + 2 \left|\beta_{i-1}\right| \theta \pi_{i-1} + \omega_i + 2 \left|\beta_{i-1}\right| \sigma_{i-1}$
		\State $e_{i-1}^j := (\mu \sqrt{n} + 2) \theta \nu_i + 2 \left|\beta_{i-1}\right| \theta \phi_{i-1} + 2 \left|\beta_{i-1}\right| \psi_{i-1}$
		\EndIf
		\If{$i=1$ \textbf{or} replace := true}
		\State $f_{i} := \epsilon \sqrt{(\mu \sqrt{n} + 1) \theta \chi_i + \zeta} + \epsilon \sqrt{ \left|\alpha_{i-1}\right| \mu \sqrt{n} \theta \pi_i} + \sqrt e_{i-1}^f \epsilon$
		\State $g_{i-1} := \epsilon \sqrt{\mu \sqrt{n} \theta \pi_i}$
		\State $h_{i} := \epsilon \sqrt{\mu \sqrt{n} \theta \xi_i} + \epsilon \sqrt{ \left|\alpha_{i-1}\right| \mu \sqrt{n} \theta \phi_i} + \sqrt e_{i-1}^h \epsilon$
		\State $j_{i-1} := \epsilon \sqrt{\mu \sqrt{n} \theta \phi_i}$
		\State replace := false
		\Else
		\State $f_{i} := f_{i-1} + \left|\alpha_{i-1}\right| \left|\beta_{i-1}\right| g_{i-2} + \left|\alpha_{i-1}\right| h_{i-1} + \sqrt e_{i-1}^f \epsilon + \left|\alpha_{i-1}\right| \sqrt e_{i-1}^g \epsilon$
		\State $g_{i-1} :=  \left|\beta_{i-1}\right| g_{i-2} + h_{i-1} + \sqrt e_{i-1}^g \epsilon$
		\State $h_{i} := h_{i-1} + \left|\alpha_{i-1}\right| \left|\beta_{i-1}\right| j_{i-2} + \sqrt e_{i-1}^h \epsilon + \left|\alpha_{i-1}\right| \sqrt e_{i-1}^j \epsilon$
		\State $j_{i-1} :=  \left|\beta_{i-1}\right| j_{i-2} + \sqrt e_{i-1}^j \epsilon$
		\EndIf
		\If{$f_{i-1} \leq \tau \rho_{i}$ \textbf{and} $f_i > \tau \rho_{i+1}$}
		\State $s_i := A p_i$; $q_i := M^{-1} s_i$; $z_i := A q_i$
    \State $r_{i+1} := b - A x_{i+1}$; $u_{i+1} := M^{-1} r_{i+1}$; $w_{i+1} := A u_{i+1}$
		\State replace := true
		\EndIf
		\EndIf
    \EndFor
    \EndProcedure 
		}
  \end{algorithmic}
\end{algorithm}

\section{Numerical results} \label{sec:numerical}

This section presents numerical results on a wide range of matrices to compare the behavior of the different 
CG methods and show the improved attainable accuracy using the automated residual replacement strategy. 
The numerical results in Sections \ref{sec:poisson} and \ref{sec:matrix} are based on a Matlab implementation of the different CG algorithms and their respective error estimates.
Parallel performance measurements in Section \ref{sec:parallel} result from a PETSc \cite{petsc-web-page} implementation of p-CG-rr on a distributed memory machine using the message passing paradigm.

\subsection{Poisson model problem} \label{sec:poisson}

\begin{table}[t]
\centering
\scriptsize
\begin{tabular}{| l | r | r r | r r | r r | r r r |}
\hline 
	\multicolumn{2}{|c|}{$\bar{x}_0 = 0$}	 & \multicolumn{2}{|c|}{CG}  & \multicolumn{2}{|c|}{CG-CG} & \multicolumn{2}{|c|}{p-CG} & \multicolumn{3}{|c|}{p-CG-rr}  \\
	\hline \hline
	Matrix   & $n$   		 & iter & relres  & iter  & relres & iter & relres & iter & relres & rr \\ 
	& &  & relerr  &  & relerr & & relerr & & relerr &   \\
\hline 
	lapl50  & 2,500    & 128  & 7.8e-15 & 127  & 8.1e-15 & 118  & 1.5e-12 & 125  & 9.1e-15 & 3  \\
	        &          &      & 6.4e-15 &      & 5.7e-15 &      & 1.1e-12 &      & 2.9e-14 &    \\ 
	lapl100 & 10,000   & 254  & 1.6e-14 & 256  & 1.6e-14 & 228  & 9.1e-12 & 272  & 1.2e-14 & 6  \\
	        &          &      & 1.4e-14 &      & 1.4e-14 &      & 6.5e-12 &      & 1.8e-14 &    \\
	lapl200 & 40,000   & 490  & 3.1e-14 & 487  & 3.2e-14 & 439  & 5.4e-11 & 536  & 2.5e-14 & 11 \\
	        &          &      & 3.7e-14 &      & 3.6e-14 &      & 5.3e-11 &      & 3.7e-14 &    \\
	lapl400 & 160,000  & 959  & 6.2e-14 & 958  & 6.4e-14 & 807  & 3.0e-10 & 957  & 4.6e-14 & 23 \\
	        &          &      & 1.0e-13 &      & 5.6e-14 &      & 3.4e-10 &      & 1.8e-13 &    \\
	lapl800 & 640,000  & 1883 & 1.2e-13 & 1877 & 1.3e-13 & 1524 & 1.4e-10 & 1876 & 1.1e-13 & 53 \\
	        &          &      & 2.7e-13 &      & 8.2e-13 &      & 2.0e-09 &      & 2.1e-13 &    \\
\hline \hline 
	\multicolumn{2}{|c|}{$\bar{x}_0 = \text{rand}(n,1)$}	 & \multicolumn{2}{|c|}{CG}  & \multicolumn{2}{|c|}{CG-CG} & \multicolumn{2}{|c|}{p-CG} & \multicolumn{3}{|c|}{p-CG-rr}  \\
	\hline \hline
	Matrix   & $n$   		 & iter & relres  & iter  & relres & iter & relres & iter & relres & rr \\
	& &  & relerr  &  & relerr & & relerr & & relerr &   \\
\hline 
	lapl50  & 2,500    & 232  & 9.0e-14 & 237  & 1.1e-13 & 197  & 1.3e-10 & 227  & 1.9e-14 & 6  \\
	        &          &      & 4.9e-14 &      & 6.4e-14 &      & 1.1e-10 &      & 7.4e-14 &    \\ 
	lapl100 & 10,000   & 444  & 2.9e-13 & 449  & 3.8e-13 & 367  & 1.5e-09 & 483  & 1.6e-14 & 10 \\
	        &          &      & 1.6e-13 &      & 2.1e-13 &      & 1.3e-09 &      & 1.5e-14 &    \\
	lapl200 & 40,000   & 881  & 1.3e-12 & 883  & 1.6e-12 & 685  & 1.7e-08 & 952  & 3.3e-14 & 20 \\
	        &          &      & 6.2e-13 &      & 7.8e-13 &      & 1.6e-08 &      & 3.9e-14 &    \\
	lapl400 & 160,000  & 1676 & 4.9e-12 & 1714 & 6.1e-12 & 1220 & 1.8e-07 & 1846 & 1.1e-13 & 35 \\
	        &          &      & 2.3e-12 &      & 2.9e-12 &      & 1.8e-07 &      & 1.5e-13 &    \\
	lapl800 & 640,000  & 3339 & 2.1e-11 & 3249 & 2.5e-11 & 2225 & 1.9e-06 & 3435 & 2.1e-12 & 65 \\
	        &          &      & 9.6e-12 &      & 1.2e-11 &      & 2.1e-06 &      & 2.4e-12 &    \\
\hline
\end{tabular}
\caption{Model problem 2D Laplacian operators of various sizes. A linear system with right-hand side $b = A\hat{x}$ where $\hat{x}_j = 1/\sqrt{n}$ is solved with the four presented algorithms. The initial guess is all-zero $\bar{x}_0 = 0$ (top table) and $\bar{x}_0 = \text{rand}(n,1)$ (bottom table). The number of iterations required to reach maximal attainable accuracy is given, along with the corresponding relative true residual norm $\|b-A \bar{x}_i\|_2 / \|b\|_2$ and the relative error (A-norm) $\|\hat{x} - \bar{x}_i\|_A / \|\hat{x}\|_A$. For the p-CG-rr method the number of replacement steps $rr$ is indicated.}
\label{tab:laplacian}
\end{table}

The methods presented above are tested on a two-dimensional Laplacian PDE model with homogeneous Dirichlet boundary conditions, discretized using second order finite differences on a uniform $n = n_x \times n_y$ point discretization of the unit square. The Poisson problem forms the basis for many practical HPC applications to which the pipelined CG method can be applied. Due to the very nature of the Laplace operator's spectrum, the application of CG to the Poisson problem typically does not display delayed convergence \cite{greenbaum1992predicting,strakovs2002error}, see \cite{gergelits2014composite} for more details. This allows us to focus on the issue of reduced attainable accuracy that is observed when applying p-CG to large scale Poisson problems in this test case.

Table \ref{tab:laplacian} shows convergence results for solving the discrete Poisson system for $n_x = n_y = 50, 100, 200, 400$ and $800$, with condition numbers ranging from $1.5\text{e+}3$ to $1.8\text{e+}5$. A linear system with exact solution $\hat{x}_j = 1/\sqrt{n}$ (such that $\|\hat{x}\|=1$) and right-hand side $b = A\hat{x}$ is solved for each of these discretization matrices. The initial guess is $\bar{x}_0 = 0$ (top table) and $\bar{x}_0 = \text{rand}(n,1)\sim U([0,1])$ (bottom table) respectively. The iteration is stopped when maximal attainable accuracy, i.e., $\min_i \|b-A\bar{x}_i\|_2$, is reached. This implies that a different stopping tolerance is used for each matrix. No preconditioner is used for the Laplace problems. The table lists the required number of iterations $iter$, the final relative true residual norm $relres$ and the final relative error $relerr$ for the CG, CG-CG, p-CG and p-CG-rr methods. Pipelined CG stagnates at a significantly larger residual and error compared to CG and CG-CG, see Section \ref{sec:analysis}. Note that for larger systems the loss of accuracy is dramatically more pronounced. 

\begin{figure}
\begin{center}
\begin{tabular}{cc}
\includegraphics[width=0.47\textwidth]{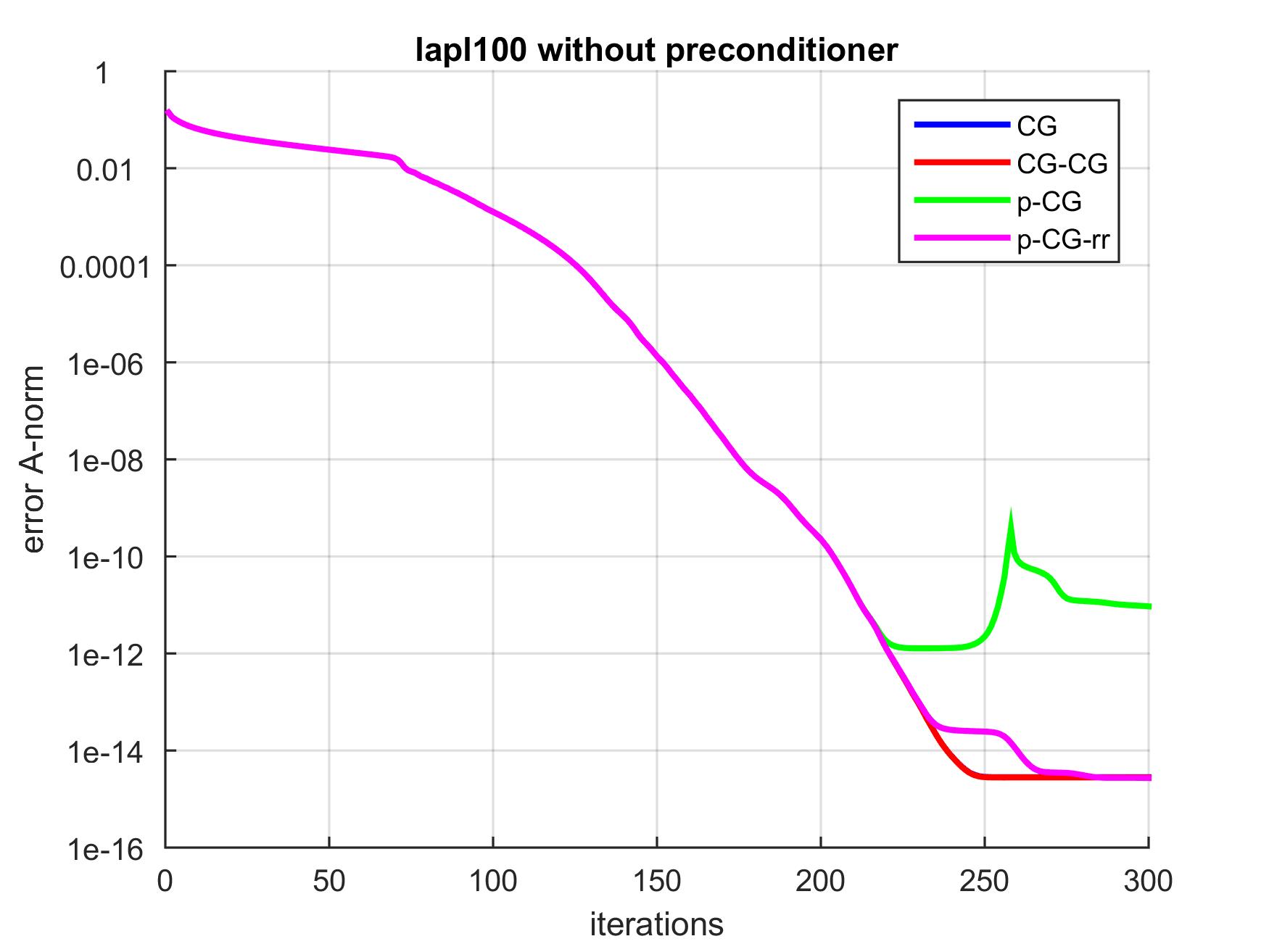} &
\includegraphics[width=0.47\textwidth]{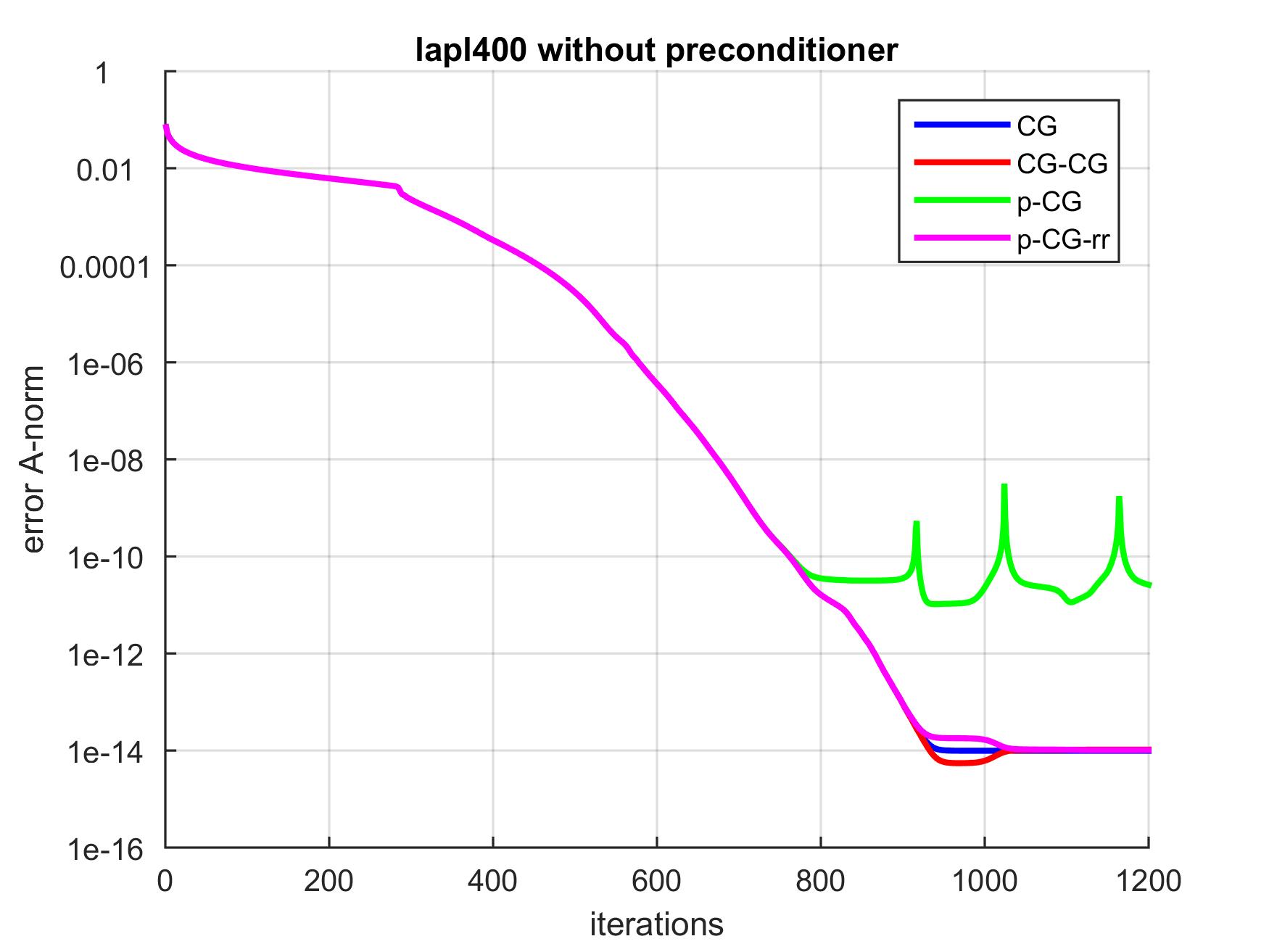} 
\end{tabular}
\end{center}
\caption{Error history for the different CG methods applied to the \texttt{lapl100} (left) and \texttt{lapl400} (right) matrices, see Table \ref{tab:laplacian}. 
Error $A$-norm $\|\hat{x} - \bar{x}_i\|_A$ as a function of iterations for CG (blue), Chronopoulos/Gear CG (red), pipelined CG (green) and p-CG-rr (magenta).
\label{fig:errAnorm_lapl}}
\end{figure}

Figure \ref{fig:errAnorm_lapl} shows the $A$-norm of the error as a function of iterations for the \texttt{lapl100} and \texttt{lapl400} problems from Table \ref{tab:laplacian}. The CG method minimizes this quantity over the respective Krylov subspace in each iteration, which (in exact arithmetic) results in a monotonically decreasing error norm. For the pipelined CG method, the error norm behaves similar to the CG method up to its stagnation point. Beyond this point the error norm is no longer guaranteed to decrease. Periodic replacement of the residual and auxiliary variables improves the attainable accuracy, as illustrated by the monotonically decreasing p-CG-rr errors. However, a slight delay of convergence \cite{strakovs2002error} is observed for the p-CG-rr method compared to classical CG, see also Table \ref{tab:laplacian}. We discuss the effect of rounding errors on orthogonality and the resulting delay of convergence near the end of Section \ref{sec:matrix}. Since the error is in general unavailable in practice, the remaining experiments focus on the norm of the residual instead.

\begin{figure}
\begin{center}
\begin{tabular}{cc}
\includegraphics[width=0.47\textwidth]{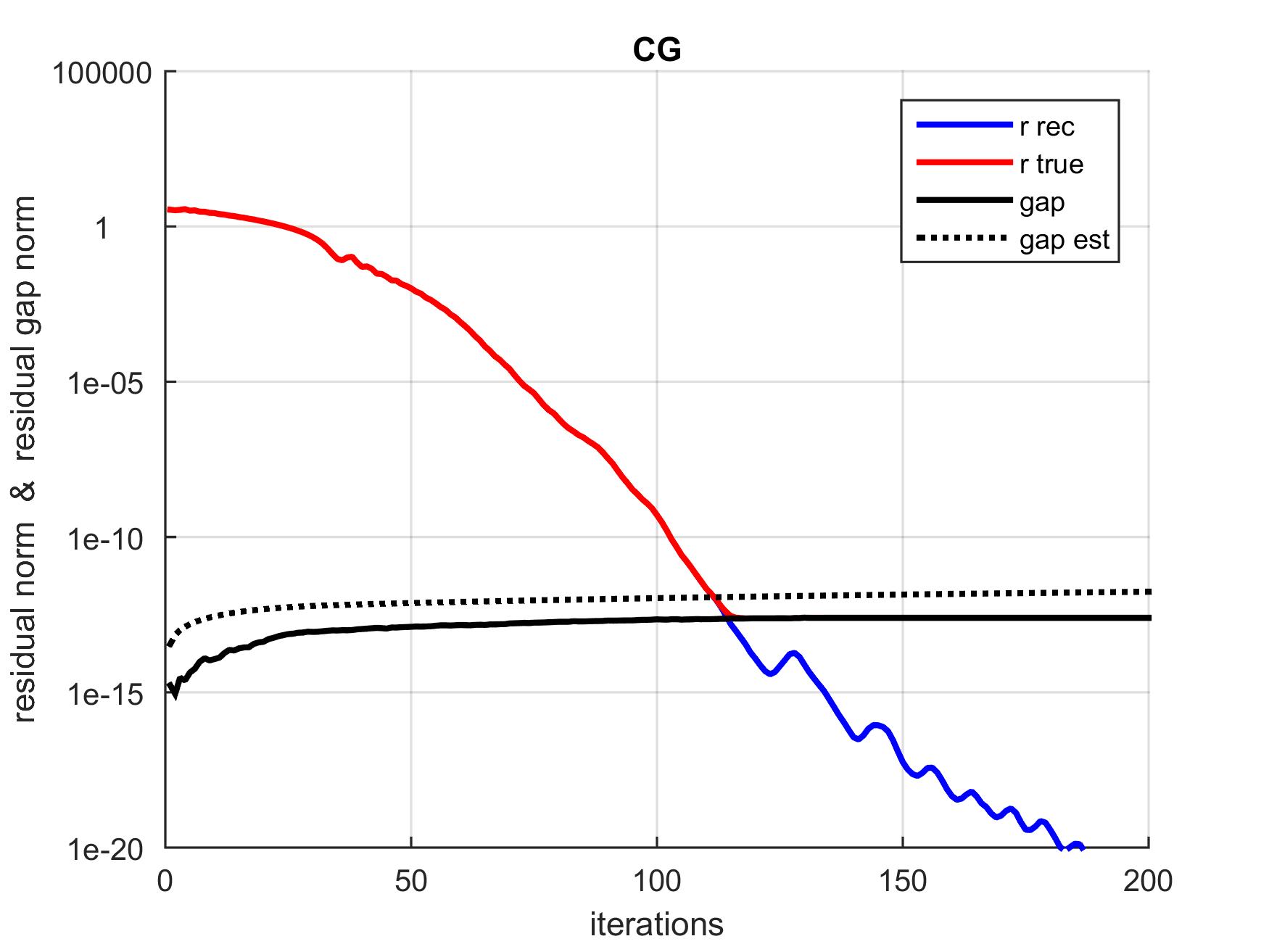} &
\includegraphics[width=0.47\textwidth]{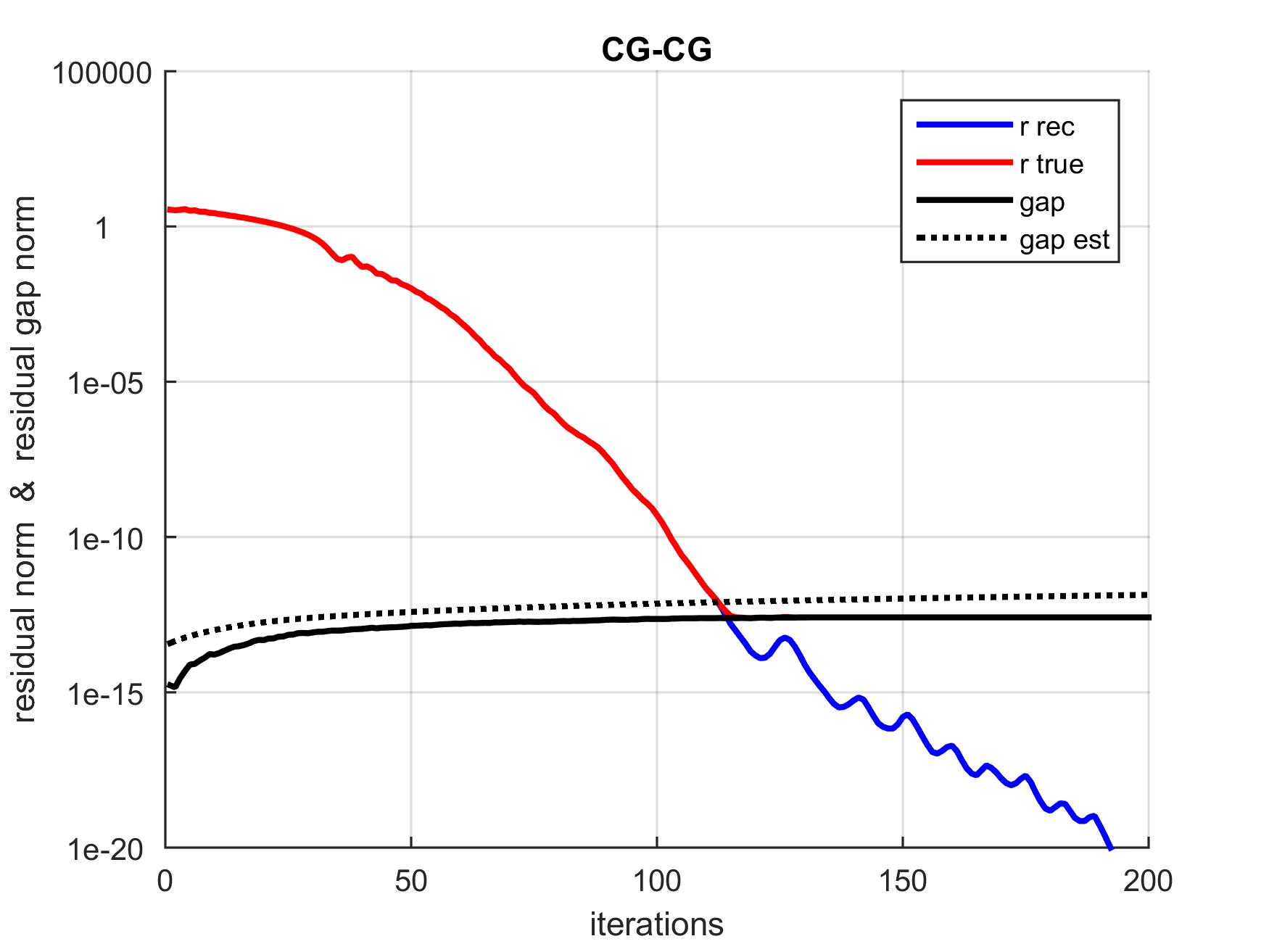} \vspace{0.5cm} \\ 
\includegraphics[width=0.47\textwidth]{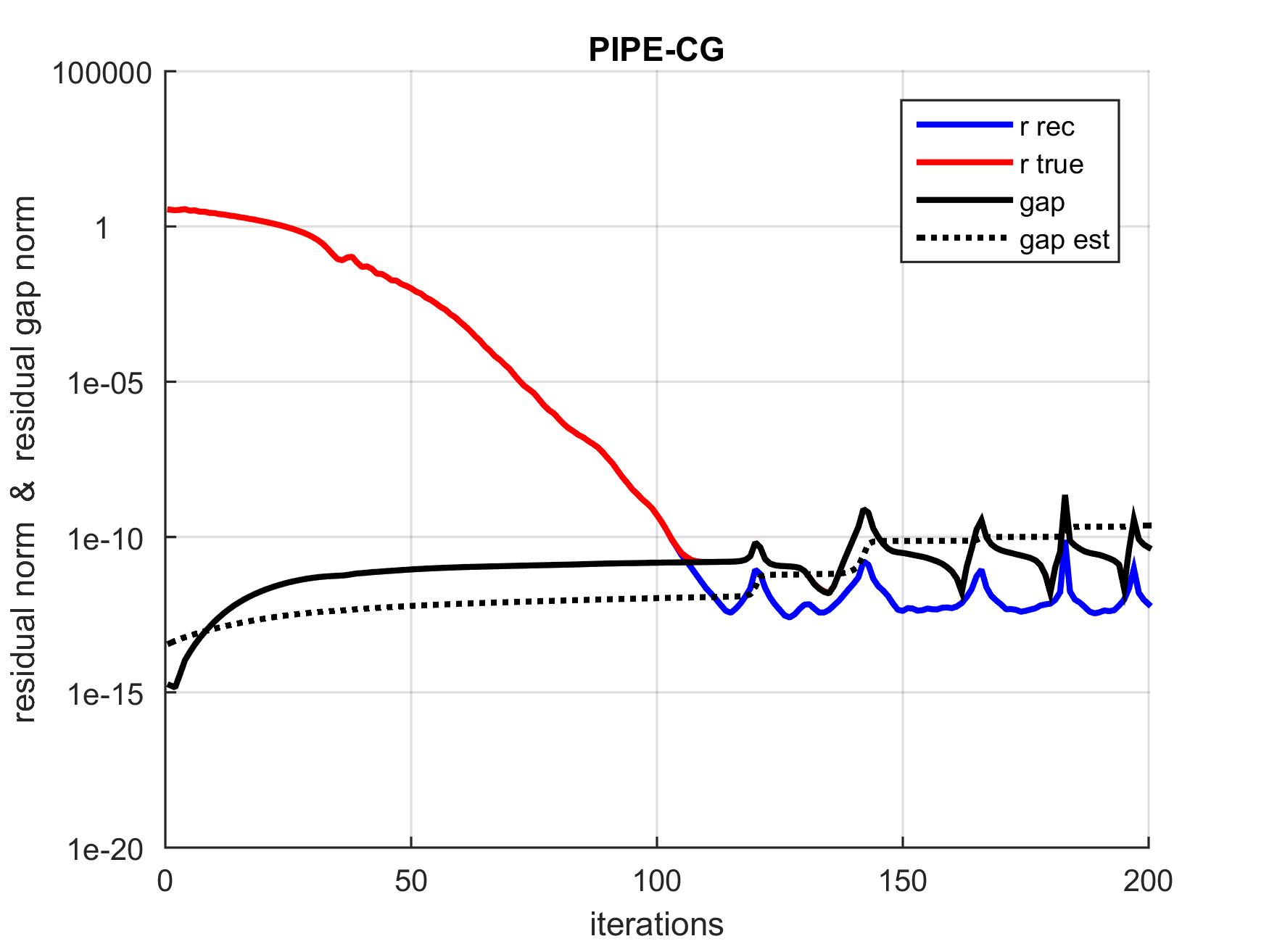} &
\includegraphics[width=0.47\textwidth]{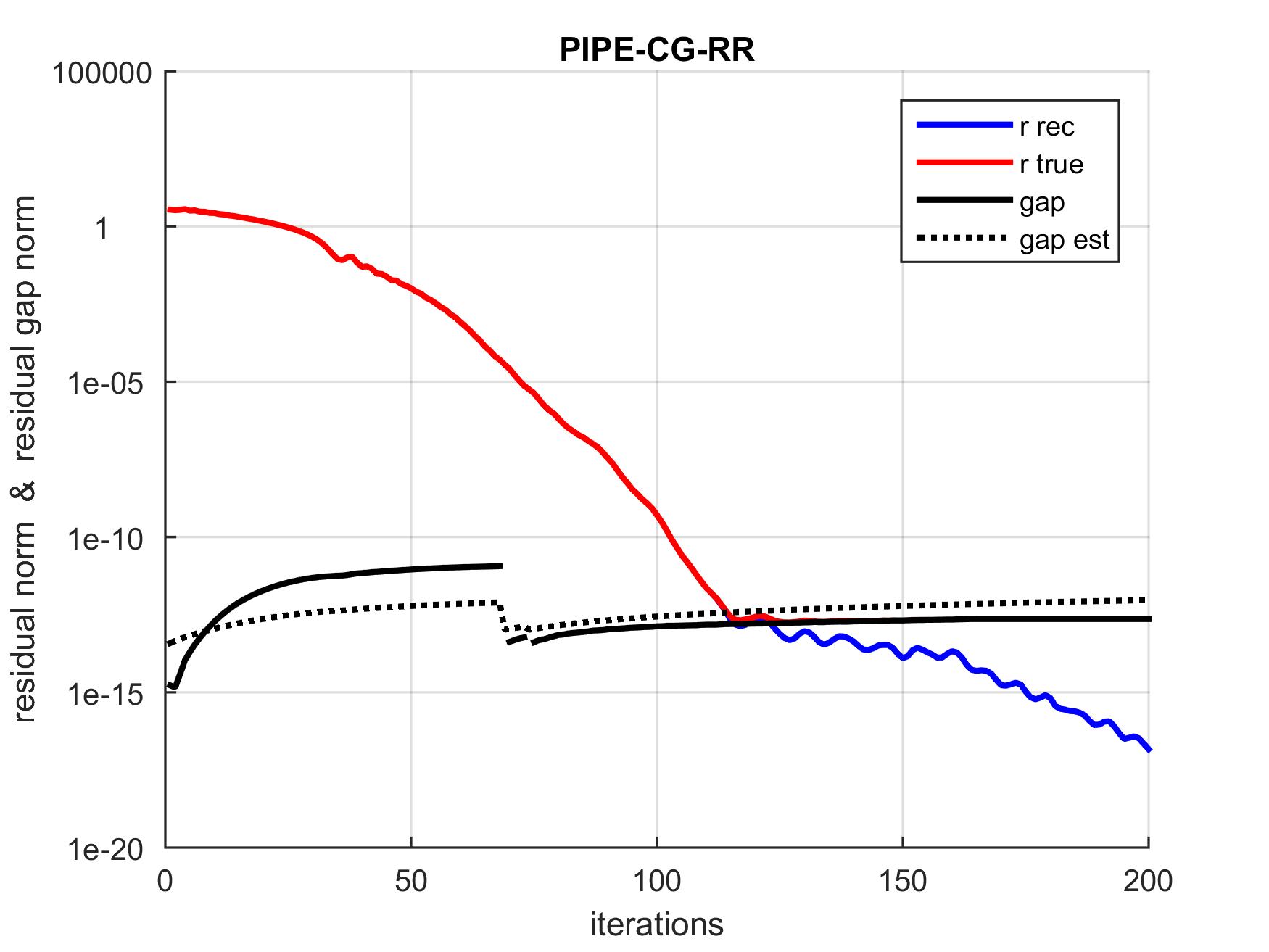}
\end{tabular}
\end{center}
\caption{Residual norm history for the different CG methods applied to the \texttt{lapl50} matrix. 
Blue: recursive residual norms $\|\bar{r}_i\|_2$. Red: true residual norms $\|b-A\bar{x}_i\|_2$. Solid black lines: residual gap norms $\|(b-A\bar{x}_i)-\bar{r}_i\|_2$ (computed explicitly). Dotted black lines: residual gap estimates $\|f_i\|_2$ computed using the expression \eqref{eq:pipecg_estimate}. 
\label{fig:conv_lapl}}
\end{figure}

Figure \ref{fig:conv_lapl} illustrates the residual convergence history and the corresponding gap between the explicit and recursive residual in the different algorithms for the \texttt{lapl50} matrix. The right-hand side is $b_j = 1/\sqrt{n}$ in this experiment. The residual norm history of CG and CG-CG is nearly indistinguishable; the norm of the true residual at stagnation is $2.4\text{e-}13$ and $2.5\text{e-}13$ respectively. The pipelined CG method suffers from the amplification of local rounding errors, leading to early stagnation of the residual norm at $1.6\text{e-}11$. The residual replacement strategy reduces the accuracy loss with a residual norm of $2.1\text{e-}13$ which is comparable to classical CG. 

\subsection{Problems from Matrix Market} \label{sec:matrix}

Numerical results on various linear systems are presented to show the effectiveness of pipelined CG with automated residual replacements. 
Table \ref{tab:matrix_market} lists all real, non-diagonal and symmetric positive definite matrices from Matrix Market\footnote{\url{http://math.nist.gov/MatrixMarket/}}, 
with their 
condition number $\kappa$, number of rows $n$ and total number of nonzero elements $\#nnz$. We solve a linear system with exact solution 
$\hat{x}_j = 1/\sqrt{n}$ and right-hand side $b = A\hat{x}$ 
with the four presented methods, using an all-zero initial guess $\bar{x}_0 = 0$.
Jacobi diagonal preconditioning (JAC) and Incomplete Cholesky Factorization (ICC) are included to reduce the number of Krylov iterations if possible. 
For the preconditioners designated by $^*$ICC an compensated Incomplete Cholesky factorization is performed, 
where a real non-negative scalar $\eta$ is used as a global diagonal shift in forming the Cholesky factor.
For the \texttt{nos1} and \texttt{nos2} matrices the shift is $\eta = 0.5$, whereas for all other $^*$ICC preconditioners we used $\eta = 0.1$.

\begin{figure}
\begin{center}
\begin{tabular}{cc}
\includegraphics[width=0.47\textwidth]{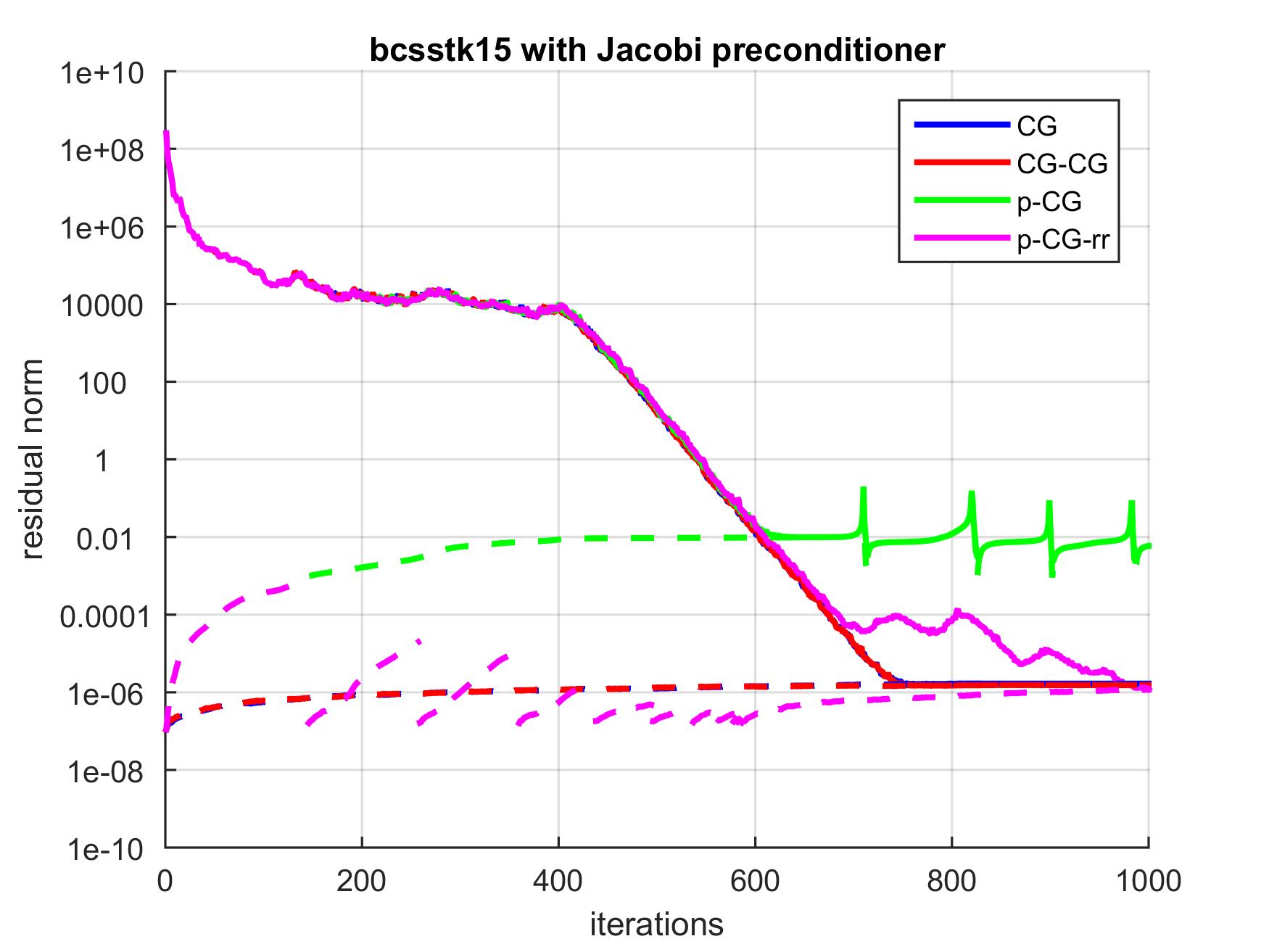} &
\includegraphics[width=0.47\textwidth]{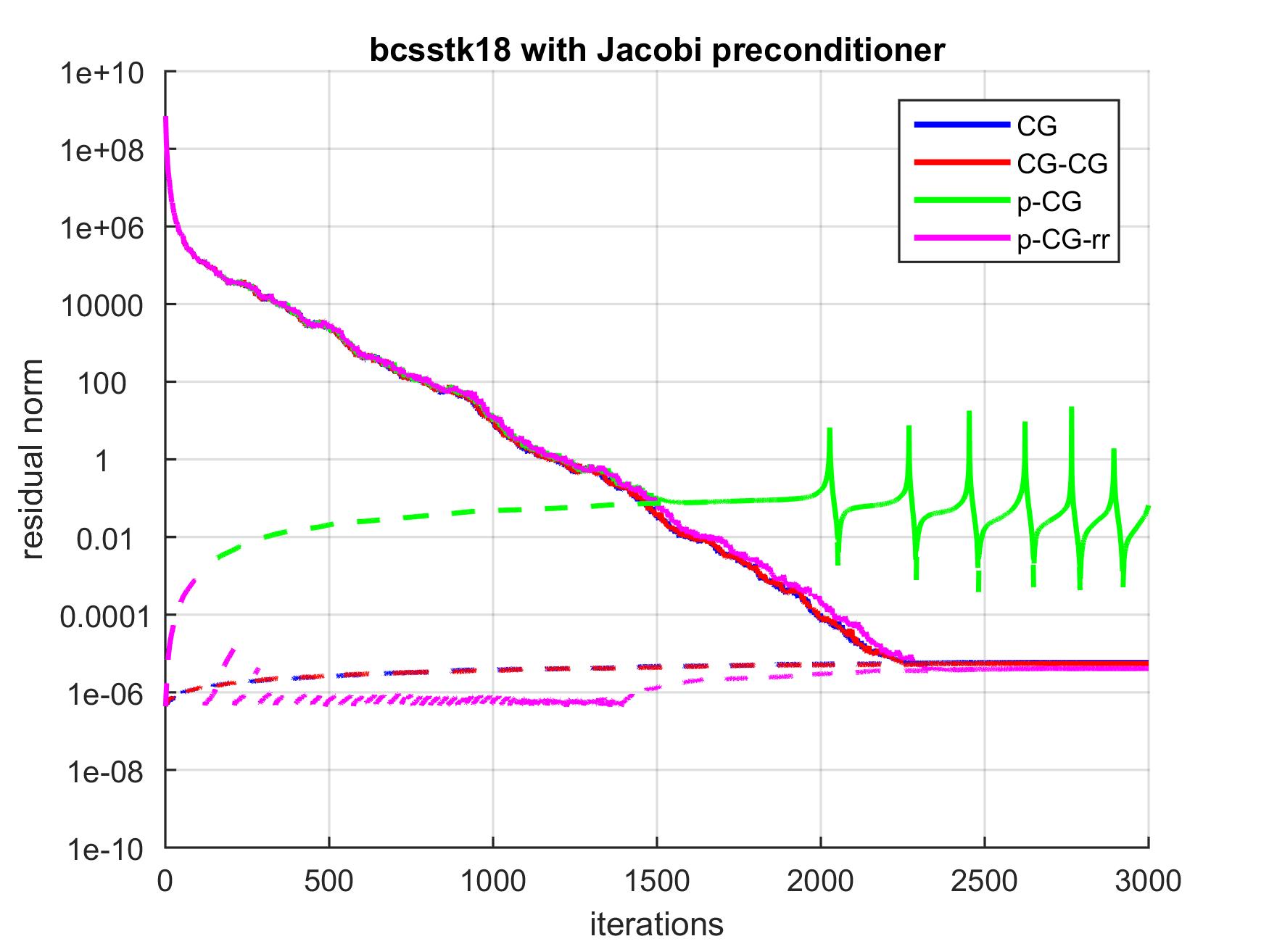} \vspace{0.5cm} \\ 
\includegraphics[width=0.47\textwidth]{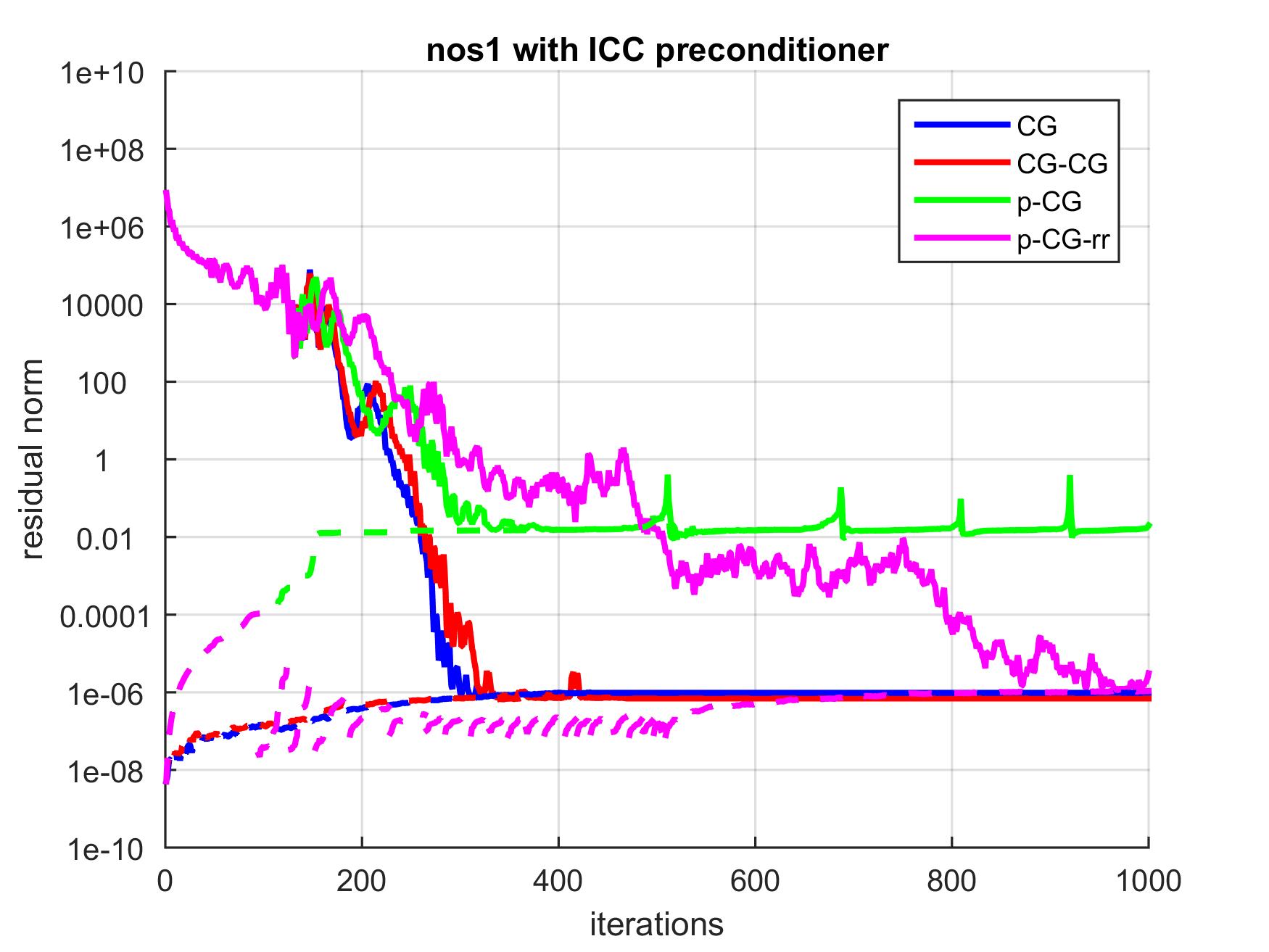} &
\includegraphics[width=0.47\textwidth]{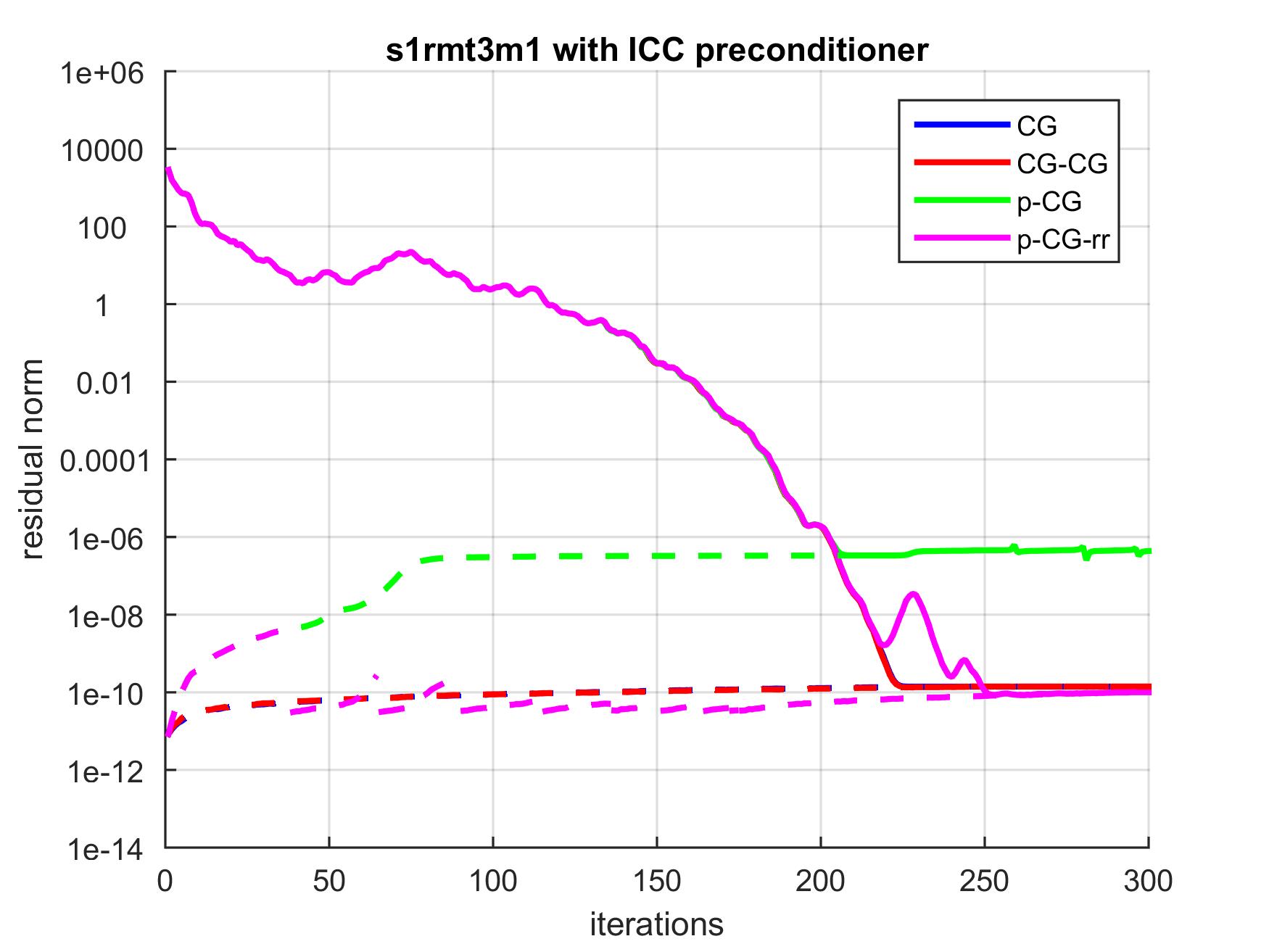}
\end{tabular}
\end{center}
\caption{Residual norm history for the different CG methods applied to four symmetric positive definite test matrices from Table \ref{tab:matrix_market}. 
Solid lines: true residual norm $\|b-A\bar{x}_i\|_2$; dashed lines: residual gap $\|f_i\|_2$. 
Convergence of CG (blue) and Chronopoulos/Gear CG (red) is largely comparable. The pipelined CG method 
(green) suffers from rounding error propagation. Automated residual replacement (magenta) 
reduces the rounding errors, leading to an accuracy that is comparable to classical CG. 
\label{fig:conv_examples}}
\end{figure}

\begin{sidewaystable}
\centering
\vspace{1.0cm}
\scriptsize
\begin{tabular}{| l | r r r r | r | r r | r r | r r | r r r |}

\hline 
	Matrix 	 & Prec  		& $\kappa(A)$ & $n$ & \#$nnz$  & $\|b\|_2$ & \multicolumn{2}{|c|}{CG}  & \multicolumn{2}{|c|}{CG-CG} & \multicolumn{2}{|c|}{p-CG} & \multicolumn{3}{|c|}{p-CG-rr}  \\
	       	 &  		 		&             &     &      	   &               & iter & relres  & iter  & relres & iter & relres & iter 		& relres & rr \\
\hline \hline
	bcsstk14 & JAC 			& 1.3e+10 & 1806   & 63,454    & 2.1e+09       & 650  & 7.6e-16 & 658  & 7.1e-16 & 506  & 5.2e-12 & 658  & 5.2e-16 & 9  \\
	bcsstk15 & JAC 			& 8.0e+09 & 3948   & 117,816   & 4.3e+08       & 772  & 3.7e-15 & 785  & 3.5e-15 & 646  & 2.3e-11 & 974  & 4.0e-15 & 10  \\
	bcsstk16 & JAC 			& 65      & 4884   & 290,378   & 1.5e+08       & 298  & 3.5e-15 & 300  & 4.0e-15 & 261  & 8.7e-12 & 301  & 2.1e-15 & 4  \\
	bcsstk17 & JAC 			& 65      & 10,974 & 428,650   & 9.0e+07       & 3547 & 1.0e-14 & 3428 & 1.7e-14 & 2913 & 2.8e-09 & 4508 & 1.2e-14 & 54 \\
	bcsstk18 & JAC 			& 65      & 11,948 & 149,090   & 2.6e+09       & 2299 & 2.2e-15 & 2294 & 2.1e-15 & 1590 & 2.9e-11 & 2400 & 1.5e-15 & 50 \\
	bcsstk27 & JAC 			& 7.7e+04 & 1224   & 56,126    & 1.1e+05       & 345  & 3.2e-15 & 345  & 4.0e-15 & 295  & 8.0e-12 & 342  & 2.7e-15 & 6  \\
	gr\_30\_30 &  -  		& 3.8e+02 & 900    & 7744      & 1.1e+00       & 56   & 2.7e-15 & 55   & 3.1e-15 & 52   & 2.0e-13 & 61   & 3.0e-15 & 2  \\
	nos1     & *ICC     & 2.5e+07 & 237    & 1017      & 5.7e+07       & 301  & 1.3e-14 & 338  & 1.2e-14 & 337  & 2.6e-10 & 968  & 1.9e-14 & 21  \\
	nos2     & *ICC     & 6.3e+09 & 957    & 4137      & 1.8e+09       & 3180 & 8.3e-14 & 3292 & 1.1e-13 & 2656 & 1.2e-07 & 4429 & 2.7e-11 & 113 \\
	nos3     & ICC 			& 7.3e+04 & 960    & 15,844    & 1.0e+01       & 64   & 1.0e-14 & 63   & 1.1e-14 & 59   & 1.0e-12 & 61   & 2.5e-14 & 3  \\
	nos4     & ICC 			& 2.7e+03 & 100    & 594       & 5.2e-02       & 31   & 1.9e-15 & 31   & 1.9e-15 & 29   & 4.0e-14 & 33   & 1.3e-15 & 2  \\
	nos5     & ICC 			& 2.9e+04 & 468    & 5172      & 2.8e+05       & 63   & 3.2e-16 & 64   & 3.4e-16 & 58   & 4.3e-14 & 65   & 2.3e-16 & 2  \\
	nos6     & ICC 			& 8.0e+06 & 675    & 3255      & 8.6e+04       & 34   & 5.1e-15 & 35   & 6.2e-15 & 31   & 5.5e-11 & 33   & 1.0e-14 & 2  \\
	nos7     & ICC 			& 4.1e+09 & 729    & 4617      & 8.6e+03       & 29   & 4.0e-14 & 31   & 2.8e-14 & 29   & 4.5e-14 & 29   & 3.0e-14 & 3  \\
	s1rmq4m1 & ICC 			& 1.8e+06 & 5489   & 262,411   & 1.5e+04       & 122  & 4.3e-15 & 122  & 4.6e-15 & 114  & 5.5e-12 & 135  & 3.7e-15 & 6  \\
	s1rmt3m1 & ICC 			& 2.5e+06 & 5489   & 217,651   & 1.5e+04       & 229  & 9.3e-15 & 228  & 8.7e-15 & 213  & 2.2e-11 & 240  & 1.7e-14 & 9  \\
	s2rmq4m1 & *ICC     & 1.8e+08 & 5489   & 263,351   & 1.5e+03       & 370  & 6.7e-15 & 387  & 7.3e-15 & 333  & 2.7e-10 & 349  & 2.5e-13 & 25  \\
	s2rmt3m1 & ICC 			& 2.5e+08 & 5489   & 217,681   & 1.5e+03       & 285  & 8.7e-15 & 283  & 1.0e-14 & 250  & 7.3e-10 & 425  & 8.7e-15 & 17  \\
	s3dkq4m2 & *ICC	    & 1.9e+11 & 90,449 & 2,455,670 & 6.8e+01       & -    & 1.9e-08 & -    & 2.1e-08 & -    & 2.8e-07 & -    & 5.6e-08 & 199 \\
	s3dkt3m2 & *ICC	    & 3.6e+11 & 90,449 & 1,921,955 & 6.8e+01       & -    & 2.9e-07 & -    & 2.9e-07 & -    & 3.5e-07 & -    & 2.9e-07 & 252 \\
	s3rmq4m1 & *ICC	    & 1.8e+10 & 5489   & 262,943   & 1.5e+02       & 1651 & 1.5e-14 & 1789 & 1.6e-14 & 1716 & 2.6e-08 & 1602 & 5.3e-10 & 154 \\
	s3rmt3m1 & *ICC	    & 2.5e+10 & 5489   & 217,669   & 1.5e+02       & 2282 & 2.7e-14 & 2559 & 2.9e-14 & 2709 & 9.3e-08 & 3448 & 8.0e-10 & 149 \\
	s3rmt3m3 & *ICC	    & 2.4e+10 & 5357   & 207,123   & 1.3e+02       & 2862 & 3.3e-14 & 2798 & 3.4e-14 & 3378 & 2.0e-07 & 2556 & 7.1e-11 & 248 \\
\hline
\end{tabular}
\caption{All real, non-diagonal and symmetric positive definite matrices from Matrix Market, listed with their respective condition number $\kappa(A)$, number of rows/columns $n$ and total number of nonzeros \#$nnz$. A linear system with right-hand side $b = A\hat{x}$ where $\hat{x}_i = 1/\sqrt{n}$ is solved with each matrix with the four presented algorithms. The initial guess is all-zero $\bar{x}_0 = 0$. Jacobi (JAC) and Incomplete Cholesky (ICC) preconditioners are included where needed. The number of iterations $iter$ required to reach maximal attainable accuracy (stagnation point) and the corresponding relative true residuals $relres =\|b-A \bar{x}_i\|_2/\|b\|_2$ are shown. For the p-CG-rr method the number of replacement steps is indicated as $rr$.}
\label{tab:matrix_market}
\end{sidewaystable}

Table \ref{tab:matrix_market} lists the number of iterations $iter$ required to reach maximal accuracy and the corresponding explicitly computed relative residual norm $relres$ for all methods. A \mbox{`-'} entry denotes failure to reach maximal accuracy within 5,000 iterations, in which case the relative residual after 5,000 iterations is displayed. The table indicates that for all test problems the residual replacement strategy incorporated in p-CG-rr improves the attainable accuracy of the p-CG method. For most matrices in the table the attainable accuracy is restored to the precision achieved by the classical CG method, although in some cases the increase in attainable accuracy is less pronounced. 

Figure \ref{fig:conv_examples} illustrates the residual norm history for a few selected matrices from Table \ref{tab:matrix_market}. 
The top panels show the true residual (solid) and residual gap (dashed) for the \texttt{bcsstk15} and \texttt{bcsstk18} matrices with Jacobi preconditioning. 
The bottom panels show the residual norm history for the \texttt{nos1} and \texttt{s1rmt3m1} matrices with ICC preconditioner. 
The residuals of the p-CG method level off sooner compared to classical CG and CG-CG. 
Based on the estimated residual gap, see \eqref{eq:pipecg_estimate}, the residual replacement strategy explicitly computes the residual 
in the iterations where the criterion \eqref{eq:criterion} is satisfied, leading to a more accurate final solution. 

Note that the behavior of the p-CG and p-CG-rr residuals near the stagnation point differs slightly from the classical CG residuals.
Furthermore, we point out that for the \texttt{nos1} matrix, see Fig.~\ref{fig:conv_examples} (bottom left), as well as several other matrices from Tables \ref{tab:laplacian} and \ref{tab:matrix_market}, the p-CG and p-CG-rr methods show significantly delayed convergence. Indeed, apart from the loss of attainable accuracy, the propagation of local rounding errors in multi-term recurrence variants of iterative schemes can cause a convergence slowdown. We refer to \cite[Section 5]{strakovs2002error} and the references therein, in particular \cite{greenbaum1989behavior}, \cite{greenbaum1992predicting} and \cite{notay1993convergence}, for a more detailed discussion on delay of convergence by loss of orthogonality due to round-off. This delay translates into a larger number of iterations required to reach a certain accuracy, see Table \ref{tab:matrix_market}. Although the residual replacement strategy ensures that a high accuracy can be obtained, it does not resolve the delay of convergence, as illustrated by the numerical results in this section. Hence, when application demands, a high accuracy can always be obtained using the p-CG-rr method, but this may come at the cost of additional iterations, inducing a trade-off between accuracy and computational effort. 

\subsection{Parallel performance} \label{sec:parallel}

\begin{figure}
\begin{center}
\begin{tabular}{cc}
\includegraphics[width=0.48\textwidth]{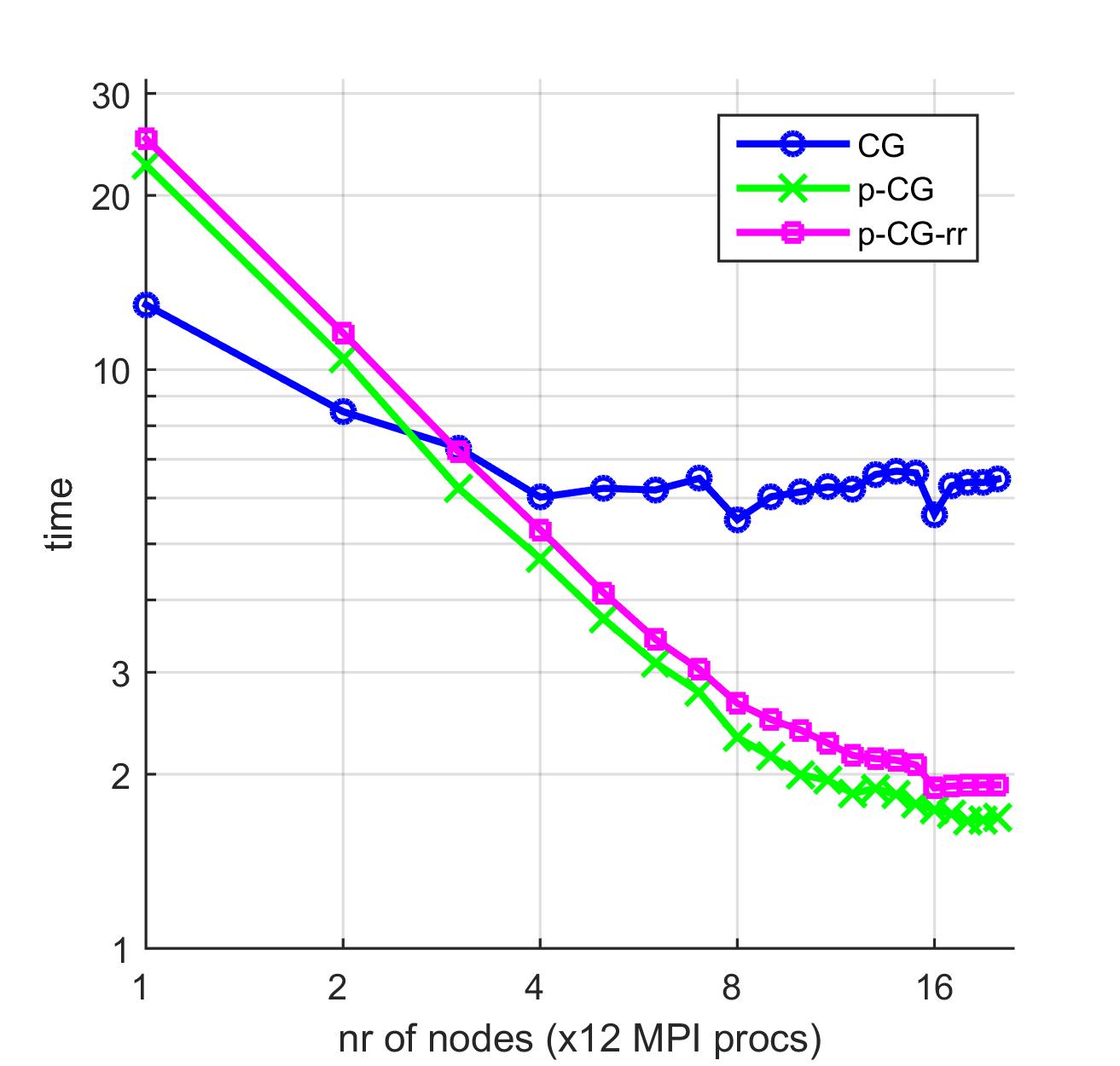} &
\includegraphics[width=0.48\textwidth]{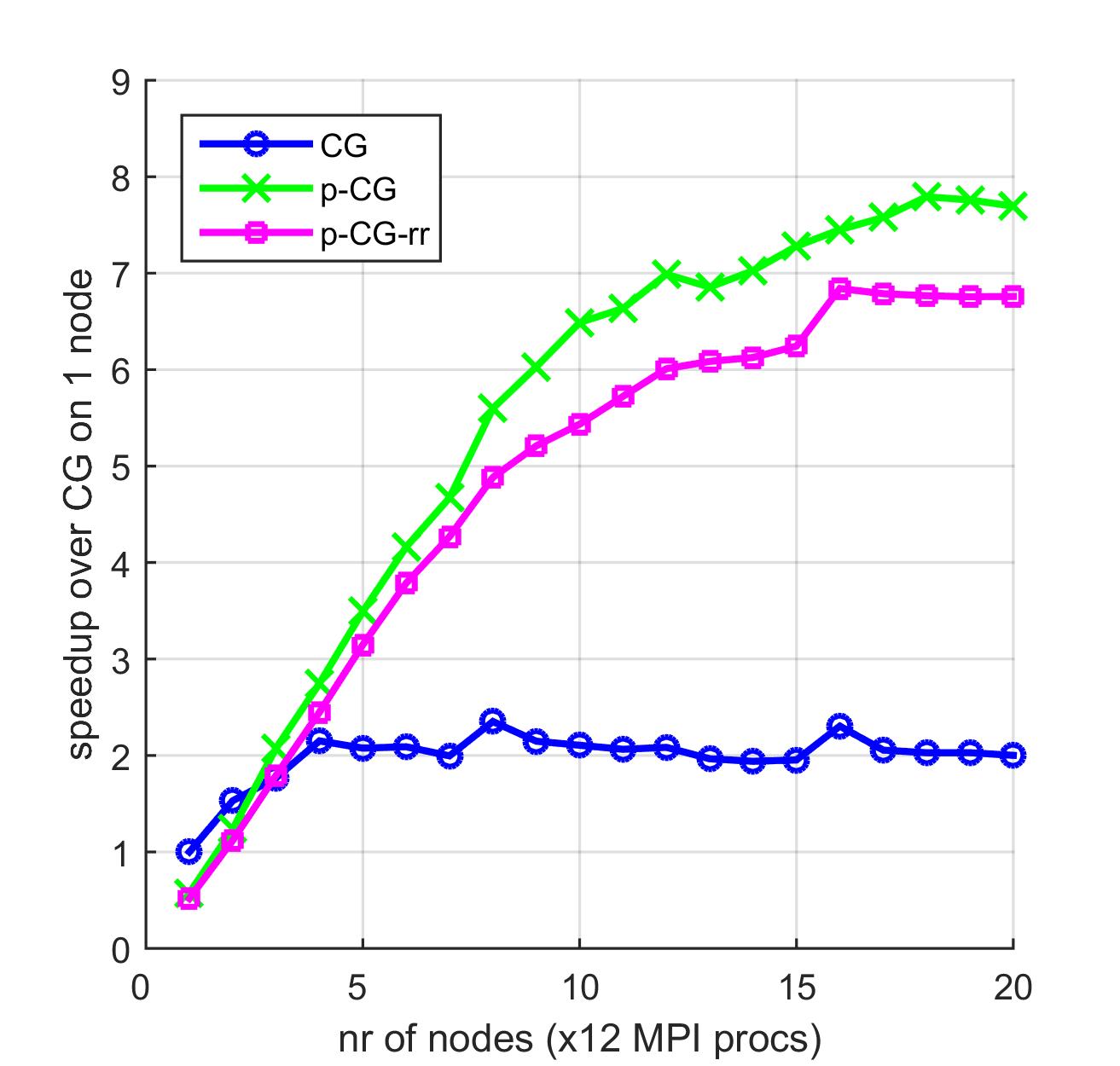} 
\end{tabular}
\end{center}
\caption{Strong scaling experiment on up to $20$ nodes ($240$ cores) for a 2D Poisson problem with $1.000.000$ unknowns.
Left: Absolute time to solution (\texttt{log10} scale) as a function of the number of nodes (\texttt{log2} scale). 
Right: Speedup over single-node classical CG. 
All methods converged in $1474$ iterations to a relative residual tolerance $1\text{e-}6$; p-CG-rr performed $39$ replacements. 
\label{fig:timings}}
\end{figure}

This section demonstrates that the parallel scalability of the pipelined CG method is maintained 
by the addition of the residual replacement strategy, and a significant speedup over classical CG can thus be obtained.
The following parallel experiments are performed on a small cluster with $28$ compute nodes, consisting of two $6$-core Intel Xeon 
X5660 Nehalem $2.80$ GHz processors each (12 cores per node). Nodes are connected by $4\,\times\,$QDR 
InfiniBand technology, providing 32 Gb/s of point-to-point bandwidth for message passing and I/O.

We use PETSc version 3.6.3. The benchmark problem used 
to asses strong scaling parallel performance is a moderately-sized 2D Poisson model, available in the PETSc 
distribution as example $2$ in the Krylov subspace solvers (KSP) folder. The simulation domain is discretized using a second order finite difference 
stencil with $1000\times1000$ grid points (1 million unknowns). No preconditioner is applied. 
The tolerance imposed on the scaled recursive residual norm $\|\bar{r}_i\|_2 / \|b\|_2$ is $10^{-6}$. 
Since each node consists of $12$ cores, we use $12$ MPI processes per node to fully exploit parallelism on the machine.
The MPI library used for this experiment is MPICH-3.1.3\footnote{\url{http://www.mpich.org/}}. Note that the environment variables 
\texttt{MPICH\_ASYNC\_PROGRESS=1} and \texttt{MPICH\_MAX\_THREAD\_SAFETY=multiple} are set to ensure optimal parallelism by allowing for non-blocking global communication.

\begin{figure}
\begin{center}
\begin{tabular}{cc}
\includegraphics[width=0.48\textwidth]{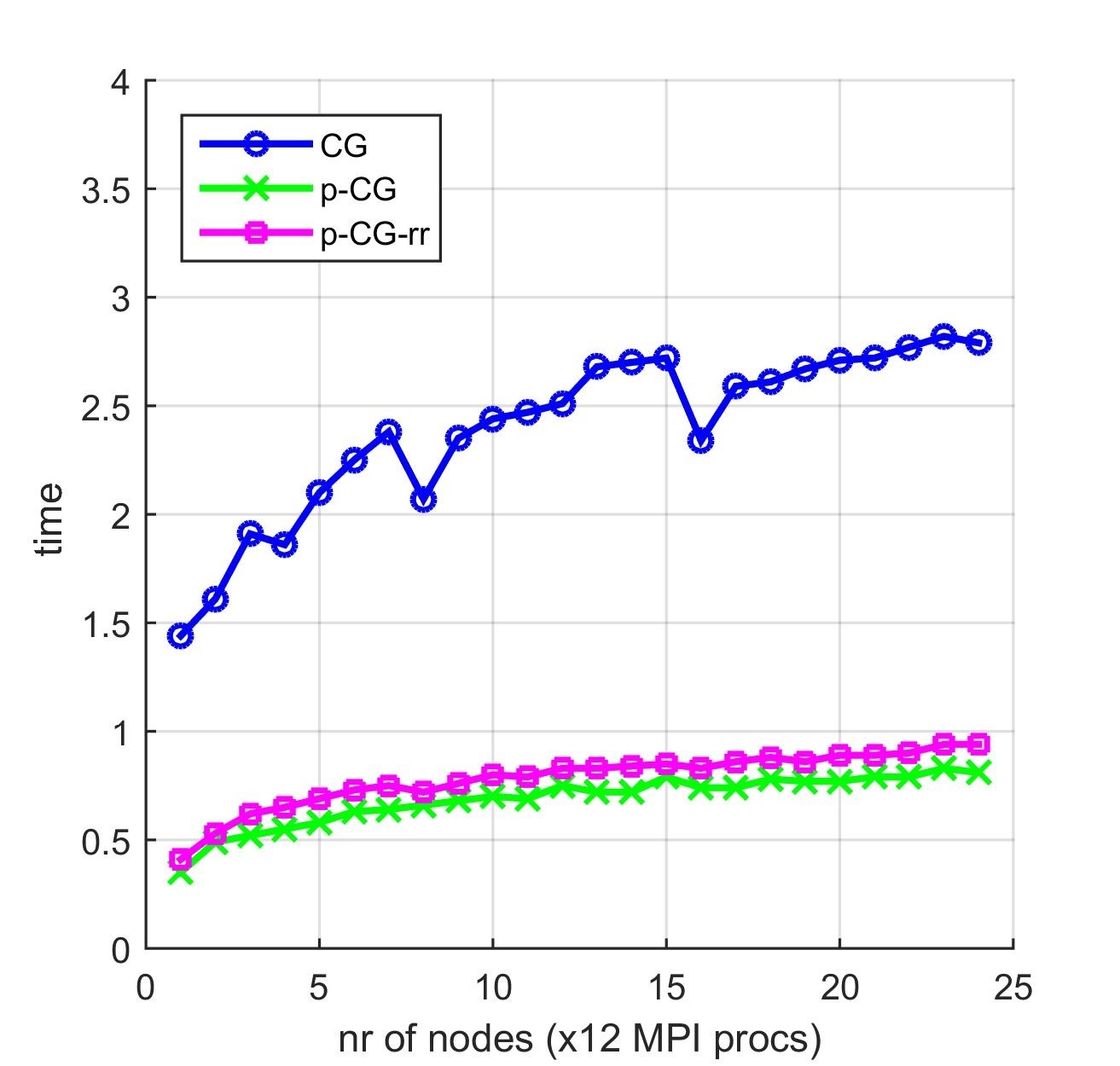} &
\includegraphics[width=0.48\textwidth]{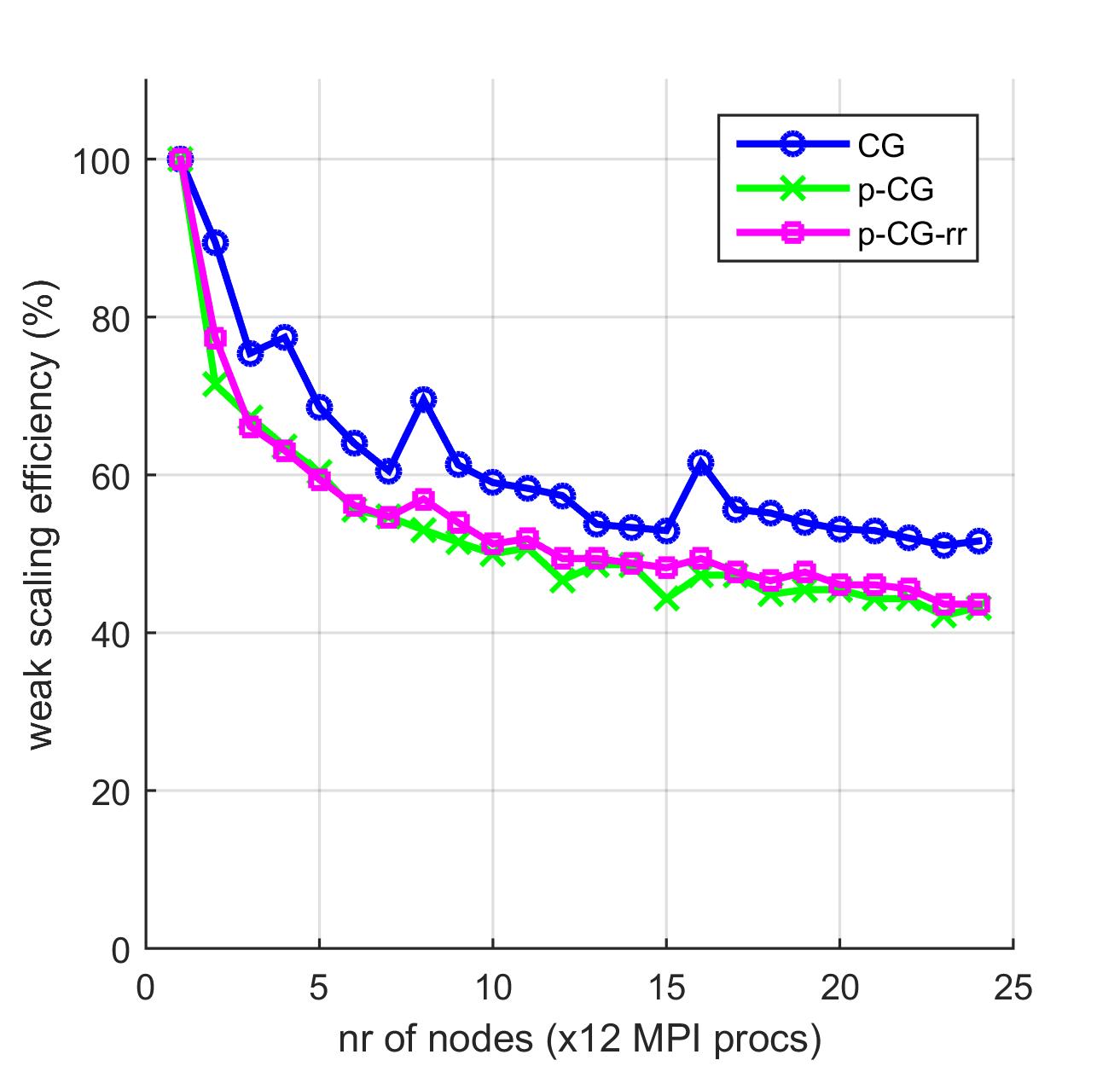} 
\end{tabular}
\end{center}
\caption{Weak scaling experiment on up to $24$ nodes ($288$ cores) for a 2D Poisson problem with $62.500$ unknowns per node ($5200$ unknowns/core).
Left: Absolute time $t_\text{CG}$ ($600$ iterations) as a function of the number of nodes. 
Right: Weak scaling efficiency relative to single-node execution: $\text{eff}_{CG}(m) = \text{t}_\text{CG}(1\,\text{node}) / \text{t}_\text{CG}(m\,\text{nodes})$. 
p-CG-rr performed $10$ replacements. 
\label{fig:timings3}}
\end{figure}

Figure \ref{fig:timings} (left) shows the time to solution as a function of the number of nodes (strong scaling). 
In this benchmark problem, pipelined CG starts to outperform classical CG when the number of nodes exceeds two. 
Classical CG stops scaling from $4$ nodes onward due to communication overhead. The pipelined methods scale well on up to 20 nodes for this problem, see Fig.~\ref{fig:timings} (right).
The maximum speedup for p-CG on $20$ nodes compared to CG on a single node is $7.7 \times$, 
whereas the CG method achieves a speedup of only $2.0 \times$ on $20$ nodes.
This implies pipelined CG attains a net speedup of $3.8 \times$ over classical CG for the current benchmark problems when both are 
executed on 20 nodes.\footnote{The theoretical time per iteration (tpi) of CG, Alg.\ref{algo::pcg}, is $2G + S$, where $G$ is the tpi 
spent by the global communication phase and $S$ is the tpi for the \textsc{spmv}. The tpi for p-CG, Alg.~\ref{algo::ppipe-cg}, is $\max(G,S)$, see \cite[Section 5]{ghysels2014hiding}. However, a third dot-product is computed in the PETSc implementations of CG and p-CG to compute the norm $\|\bar{r}_i\|_2$. The tpi for PETSc's CG is thus $3G + S$. For p-CG this extra dot-product is combined into the existing global reduction phase such that the tpi remains unaltered. Hence, when $G = S$ a theoretical maximal speedup factor of $4\times$ can be achieved by p-CG.}
Performance of p-CG-rr is comparable to that of p-CG. The minor observed slowdown  
is primarily due to the additional computational work required for the \textsc{spmv}s \eqref{eq:rrs} when replacement takes place. 
The p-CG-rr algorithm achieves a speedup of $3.4 \times$ over CG on $20$ nodes for this problem and hardware setup.
Note that the pipelined variants are effectively slower than classical CG on one or two nodes. This is due to the additional \textsc{axpy}s in the pipelined methods, which require a significant amount of time for smaller numbers of nodes but are negligible on large numbers of nodes due to parallelism. This observation illustrates that good parallel algorithms are not necessarily the best sequential ones.

Figure \ref{fig:timings3} displays results for a weak scaling experiment, where the size of the Poisson problem grows linearly 
with respect to the number of cores. A fixed number of $62.500$ unknowns per node is used. The problem hence consists of 
$1225 \times 1225$ (1.5 million) unknowns on 24 nodes. Fig.~\ref{fig:timings3} (left) shows the time required to perform $600$ iterations (fixed) 
of the various methods on up to $24$ nodes. The speedup observed for the pipelined methods in Fig.~\ref{fig:timings} is again visible here.
The weak scaling efficiency of the p-CG and p-CG-rr algorithms (relative to their respective single-node execution) on $24$ nodes ($43\%$) is comparable 
to that of classical CG ($51\%$), see Fig.~\ref{fig:timings3} (right). 

\begin{figure}
\begin{center}
\begin{tabular}{cc}
\includegraphics[width=0.48\textwidth]{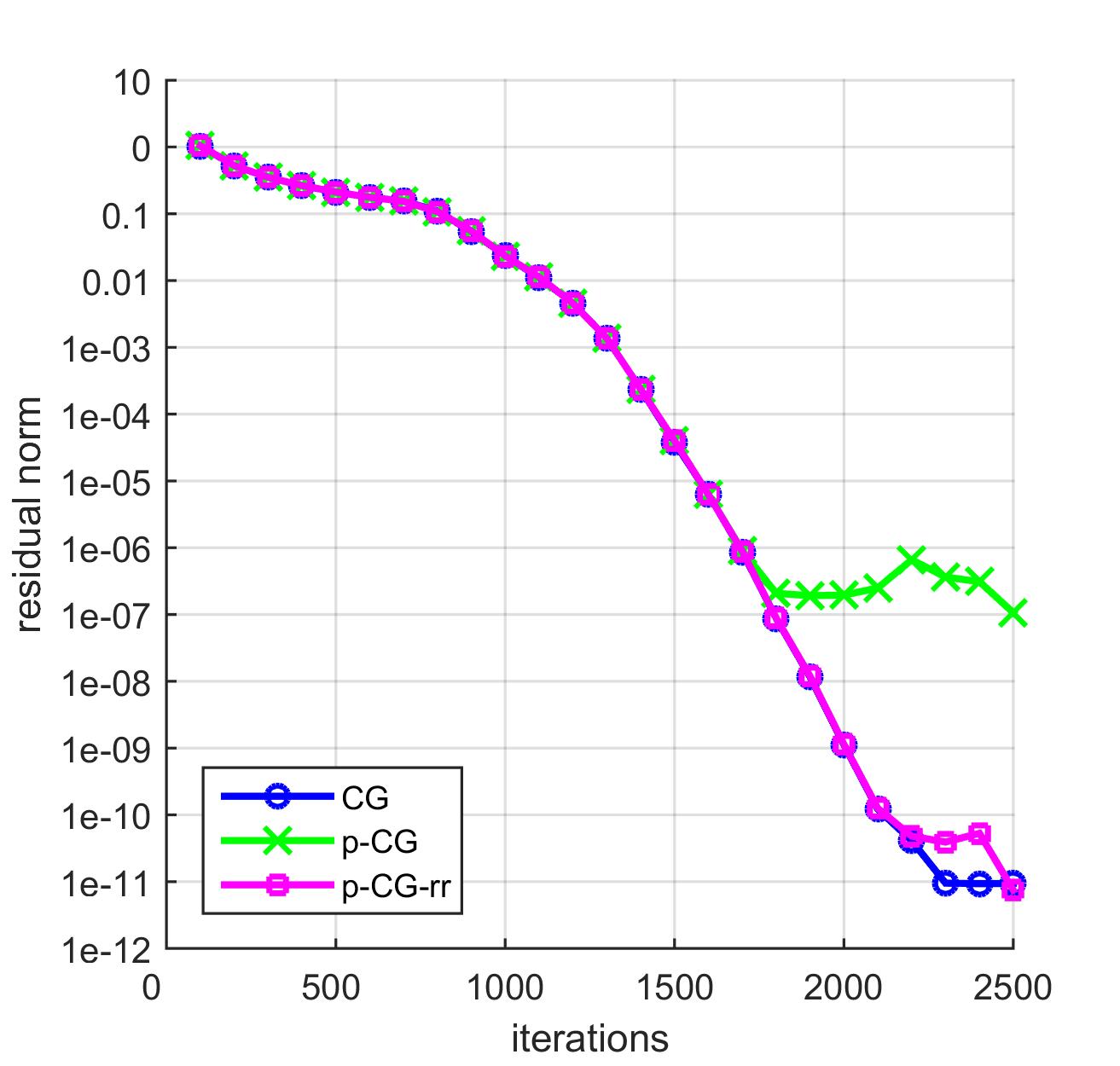} &
\includegraphics[width=0.48\textwidth]{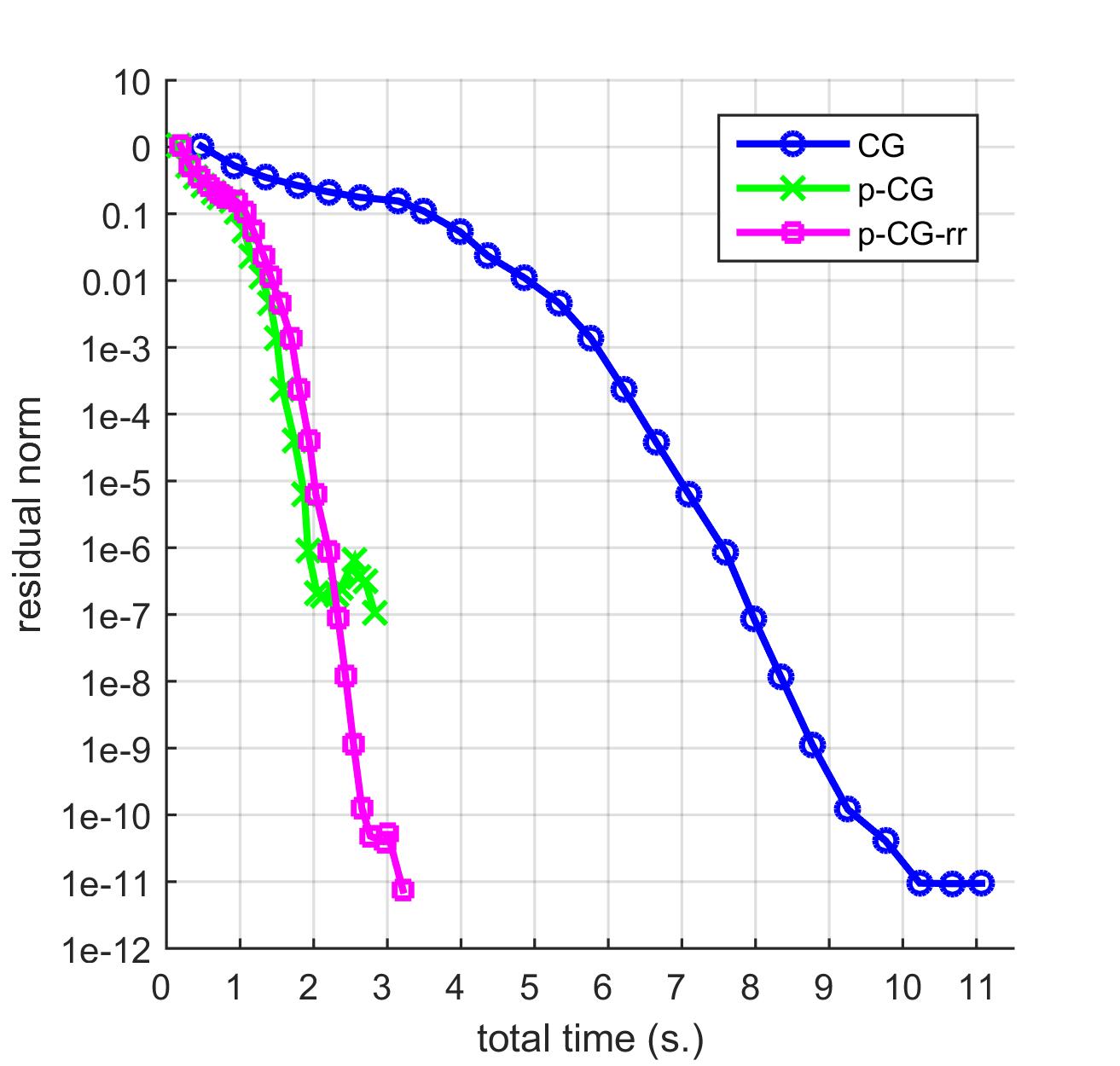} 
\end{tabular}
\end{center}
\caption{Accuracy experiment on $20$ nodes ($240$ cores) for a 2D Poisson problem with $1.000.000$ unknowns.
Left: Explicitly computed residual as a function of iterations.
Right: Residual as a function of total time spent by the algorithm. 
Maximal number of iterations is $2500$ for all methods; p-CG-rr performed (maximum) $39$ replacements.
\label{fig:timings2}}
\end{figure}

Figure \ref{fig:timings2} shows the accuracy of the solution as a function of the number of iterations (left) and computational time (right) spent
by the algorithms for the 2D Poisson $1000 \times 1000$ benchmark problem on a 20 node setup. In $3.2$ seconds ($\sim$ $2500$ iterations) 
the p-CG-rr algorithm obtains a solution with true residual norm $7.5\text{e-}12$. 
Classical CG is over three times slower, requiring 11.1 seconds to attain a comparable accuracy 
(residual norm $9.4\text{e-}12$),
see also Fig.~\ref{fig:timings}.
The p-CG method without residual replacement is unable to reach a comparable accuracy regardless of computational effort.
Indeed, stagnation of the true residual norm around $2.0\text{e-}7$ is imminent from a total time of $2.0$ seconds ($\sim$ $1800$ iterations) onward.
For completeness we note that the speedup of p-CG/p-CG-rr over classical CG can also be obtained for less accurate final solutions, 
e.g., with $\|\bar{r}_i\| = 10^{-8}$ or $10^{-6}$, as shown by Fig.~\ref{fig:timings2} (right).

\section{Conclusions} \label{sec:conclusions}

Deviation of the recursively computed residuals from the true residuals due to
the propagation of local rounding errors is a well-known issue in many numerical methods. 
This behavior is significantly more prominent in multi-term recursion variants
of CG, such as the three-term recurrence CG algorithm \cite{stiefel1955relaxationsmethoden},
Chronopoulos \& Gear's communication avoiding CG (CG-CG) \cite{chronopoulos1989s}, and 
the communication hiding pipelined CG method (p-CG) \cite{ghysels2014hiding}. For these methods, the dramatic amplification of local rounding errors may lead to a stagnation of the residual norm at several orders of magnitude above the accuracy attainable by classical CG.

This paper aims to lay the foundation for the analysis of the propagation of local rounding errors that stem from the recursions in 
classical CG, CG-CG and the pipelined CG algorithm. Pipelined CG 
features additional recursively computed auxiliary variables 
compared to classical CG which are all prone to rounding errors. We show that the gap between 
the explicitly computed and recursive residual is directly related to the gaps on the 
other recursively defined auxiliary variables. 
A bound on the residual gap norm is derived, which provides insight into the observed accuracy loss in the pipelined CG algorithm. 
Furthermore, a practically useable estimate for the residual gap is suggested. 
Based on this estimate, a heuristic to compensate for the loss of attainable accuracy is proposed in the 
form of an automated residual replacement strategy. 
However, since the assumption of Krylov basis orthogonality is not guaranteed to hold in finite precision, 
the replacement strategy should be interpreted as an effective yet primarily intuitive practical tool to improve attainable accuracy.
 

The residual replacement strategy is illustrated on a variety of numerical benchmark problems. 
For most test problems the replacements allow to attain a significantly improved accuracy which is unobtainable by pipelined CG. However, a delay of convergence due to the replacements is observed for specific problems. 
Although the incorporation of the replacement strategy in the algorithm requires 
the calculation of several additional vector norms, these computations can easily be combined 
with the existing global communication phase. 
Performance results with a parallel implementation of p-CG-rr in PETSc using the MPI message passing paradigm 
indicate that the replacement strategy improves accuracy but does not impair parallel performance.

While this work aims to be a contribution towards more efficient and more accurate parallel variants of the Conjugate Gradient algorithm, 
we point out that many open questions in this area still remain. Future work may (and should) tackle the issues of effectively accounting for the loss of orthogonality in finite precision in the numerical framework presented in this paper, as well as the analysis and possible remedy for the observed delay of convergence in multi-term recurrence CG variants, including pipelined CG.

\section{Acknowledgments} This work is funded by the EXA2CT European Project on
Exascale Algorithms and Advanced Computational Techniques, which receives 
funding from the EU's Seventh Framework Programme (FP7/2007-2013) under grant agreement no.~610741.
Additionally, S.\,C.\,is funded by the Research Foundation Flanders (FWO) under grant 12H4617N. 
The authors would like to cordially thank both Zden{\v{e}}k Strako{\v{s}} and the anonymous SIMAX referee for their useful comments and valuable suggestions on earlier versions of this manuscript.

\bibliographystyle{plain}
\bibliography{refs2}

\end{document}